# CPGD: Cadzow Plug-and-Play Gradient Descent for Generalised FRI

Matthieu Simeoni, *Member, IEEE*, Adrien Besson, *Associate Member, IEEE*, Paul Hurley, *Senior Member, IEEE*, and Martin Vetterli, *Fellow, IEEE*

***Abstract*—Finite rate of innovation (FRI) is a powerful reconstruction framework enabling the recovery of sparse Dirac streams from uniform low-pass filtered samples. An extension of this framework, called generalised FRI (genFRI), has been recently proposed for handling cases with arbitrary linear measurement models. In this context, signal reconstruction amounts to solving a joint constrained optimisation problem, yielding estimates of both the Fourier series coefficients of the Dirac stream and its so-called annihilating filter, involved in the regularisation term. This optimisation problem is however highly non convex and non linear in the data. Moreover, the proposed numerical solver is computationally intensive and without convergence guarantee. In this work, we propose an implicit formulation of the genFRI problem. To this end, we leverage a novel regularisation term which does not depend explicitly on the unknown annihilating filter yet enforces sufficient structure in the solution for stable recovery. The resulting optimisation problem is still non convex, but simpler since linear in the data and with less unknowns. We solve it by means of a provably convergent proximal gradient descent (PGD) method. Since the proximal step does not admit a simple closed-form expression, we propose an inexact PGD method, coined Cadzow plug-and-play gradient descent (CPGD). The latter approximates the proximal steps by means of Cadzow denoising, a well-known denoising algorithm in FRI. We provide local fixed-point convergence guarantees for CPGD. Through extensive numerical simulations, we demonstrate the superiority of CPGD against the state-of-the-art in the case of non uniform time samples.**

Manuscript received May 1, 2020; revised October 12, 2020 and November 17, 2020; accepted November 21, 2020. Date of publication November 27, 2020; date of current version December 23, 2020. The associate editor coordinating the review of this manuscript and approving it for publication was Prof. Shogo Muramatsu. The work of Matthieu Simeoni was supported by the Swiss National Science Foundation under Grant 200 021 181 978/1, "SESAM - Sensing and Sampling: Theory and Algorithms". *(Matthieu Simeoni, and Adrien Besson have contributed equally to this work.) (Corresponding author: Matthieu Simeoni.)*

Adrien Besson is with the E-Scopics, 13760 Saint-Cannat, France (e-mail: adribesson@gmail.com).

Martin Vetterli is with the Laboratoire de Communications Audiovisuelles (LCAV), École Polytechnique Fédérale de Lausanne (EPFL), CH-1015 Lausanne, Switzerland (e-mail: martin.vetterli@epfl.ch).

Matthieu Simeoni was with the Foundations of Cognitive Solutions group in the IBM Research Laboratory of Zurich. He is now with the Laboratoire de Communications Audiovisuelles (LCAV), École Polytechnique Fédérale de Lausanne (EPFL), CH-1015 Lausanne, Switzerland (e-mail: matthieu.simeoni@epfl.ch).

Paul Hurley is with the Centre for Research in Data Science & Mathematics, Western Sydney University (WSU), Parramatta, NSW 2150, Australia (e-mail: Paul.Hurley@westernsydney.edu.au).

## I. Introduction

Sampling theorems lie at the foundation of modern *digital signal processing* as they permit the convenient navigation between the *analogue* and *digital* worlds [1], [2]. The most famous is undoubtedly the *Shannon sampling theorem* [3], which states that *bandlimited* signals can be recovered exactly from their discrete samples for a sufficient sampling rate. This major result has had tremendous impact on the field of signal processing and by extension on many fields of natural sciences. But this unanimous celebration lead many scientists to start thinking about sampling theory exclusively in terms of bandlimitedness, which is only a *sufficient* condition for a signal to admit a discrete representation. In fact, sampling theorems can also be devised for non-bandlimited signals as long as they possess *finitely many* degrees of freedom.

This fact was brought to the attention of the signal processing community in [4], where the authors introduced the *finite rate of innovation (FRI)* framework. FRI is concerned with the sampling of *sparse* non-bandlimited signals such as the prototypical sparse signal, namely the $T$-periodic *stream of Diracs:*

$$x(t) = \sum_{k'\in\mathbb{Z}} \sum_{k=1}^{K} x_k \delta(t - t_k - Tk'), \quad \forall t \in \mathbb{R}, \tag{1}$$

with $x_k \in \mathbb{C}$ and $t_k \in [0, T[$. In the FRI framework, the sparsity is measured in terms of its *rate of innovation*, defined as the number of degrees of freedom per unit of time. For instance, the Dirac stream (1) has $2K$ degrees of freedom $\{x_k, t_k\}_{k=1,\dots,K}$ per period $T$, yielding a finite rate of innovation of $\rho = 2K/T$. Intuitively, any lossless sampling scheme for (1) must therefore have a sampling rate at least as large as the rate of innovation $\rho$ or it will be impossible to fix all degrees of freedom. The reconstruction of FRI signals has useful applications in many fields of applied signal processing such as ultra-wide band communications [5], [6], electroencephalography (EEG) [7], [8], optical coherence tomography [9], ultrasound imaging [10]–[12], radio astronomy [13], array signal processing [14], calcium imaging [15] non uniform spline approximation [16], [17] and functional magnetic resonance imaging (fMRI) [18].

Blu *et al.* described in [19] a sampling scheme achieving the second best sampling rate after the critical innovation rate,

permitting to perfectly recover the signal innovations from the knowledge of any $2K+1$ consecutive Fourier coefficients of $x$. Unfortunately, this scheme can be very sensitive to noise perturbations in the collected samples. This is because the recovery of the innovations $t_k$ relies on the resolution of a so-called *annihilating equation* which requires the Toeplitz matrix built from the Fourier coefficients to be rank deficient. While this structural constraint is guaranteed to hold in the case of noiseless recovery of Dirac streams, it can break down in the presence of noise, inevitable in practical applications.

As a remedy, Blu *et al.* proposed to *denoise* the collected samples prior to solving the annihilating equation. To this end, they leveraged the well-known *Cadzow algorithm* [20], which aims to retrieve the closest rank-deficient Toeplitz matrix to a high-dimensional embedding of the data via an alternating projection method. When upgraded with this extra denoising step, simulation results from Blu *et al.* revealed that the overall accuracy of the recovery procedure remains very good for signal-to-noise ratios (SNR) as low as 5 dB [19]. While Cadzow algorithm empirically provides accurate results after a few iterations, convergence in theory has however not been demonstrated to date, due to the non convex nature of the space of rank-deficient matrices. Condat and Hirabayashi [21] revisited Cadzow denoising as a *structured low-rank approximation (SLRA)* problem and proposed a *Douglas-Rachford splitting* algorithm to solve it [22], with higher accuracy than traditional Cadzow denoising. Unfortunately, the gain in accuracy comes at the price of significantly higher computational cost, the Douglas-Rachford splitting method requiring many more iterations to converge than Cadzow algorithm.

In addition to their somewhat heuristic nature, neither Cadzow denoising nor its upgrade can handle more general types of input measurements as considered in the *generalised FRI (genFRI)* framework introduced by Pan *et al.* in [23]. The latter extends FRI to very generic cases where the measurements are related to the unknown Fourier coefficients of signals satisfying the annihilating property by a linear map. In such configurations, both the Fourier coefficients and their corresponding annihilating filter are unknown and must be estimated from the data. Pan *et al.* proposed to perform this joint estimation task by solving a constrained optimisation problem which recovers the Fourier coefficients, required to minimise a quadratic data-fidelity term, and their corresponding annihilating filter coefficients. The annihilating equation linking the two unknowns is explicitly enforced as a constraint. This optimisation problem is highly non convex and non linear in the data. They suggested to solve it via an iterative alternating minimisation algorithm with multiple random initialisations [23]. The proposed algorithm however comes without convergence guarantees and is computationally intensive.

In this paper, we propose an implicit formulation of the genFRI problem in which only the Fourier coefficients to be annihilated are recovered. This formulation does not rely explicitly on the unknown annihilating filter but rather leverages a structured low-rank regularisation constraint based on a "Toeplitzification" linear operator, guaranteeing non-trivial solutions to the annihilating equation. The resulting optimisation problem is still non convex, but simpler to analyse and solve since it is linear in the data and with less unknowns. We solve the implicit genFRI problem via *proximal gradient descent (PGD)* [24], [25].

We first consider PGD with *exact* proximal steps which is shown to converge towards critical points of the implicit genFRI problem. The latter is however impractical since the proximal step involved at each iteration does not have a closed-form expression. We therefore consider an *inexact* PGD [26], with proximal steps approximated by means of alternating projections. In the case of injective forward matrices, the approximate proximal step is shown to reduce to Cadzow denoising. Such an approach is reminiscent of the *plug-and-play (PnP) framework* in which proximal operators involved in first-order iterative methods are replaced by generic denoisers [27]–[29]. For this reason, we name our reconstruction algorithm *Cadzow PnP Gradient Descent (CPGD)*.[1] We demonstrate that CPGD converges locally towards fixed points of the update equation for injective forward matrices. Through simulations of irregular and noisy time sampling of periodic stream of Diracs we show that CPGD is almost always more accurate and more efficient than the procedure proposed by Pan *et al.* in [23], sometimes by several orders of magnitude.

The remainder of the paper is organised as follows:

- Preliminary concepts required for the understanding of the further sections are introduced in Section II.
- Section III describes the genFRI problem and details the proposed implicit formulation. The CPGD algorithm is introduced in Section IV.
- Experiments and results are detailed in Section V and concluding remarks are given in Section VI.

Note that Appendices D to L are provided as supplementary material to this manuscript. Finally, all experiments and simulations are fully reproducible using the benchmarking routines provided in our GitHub repository [30].

## II. Preliminaries

In this section we introduce a linear operator, baptised *Toeplitzification operator*,[2] which transforms a vector into a Toeplitz matrix. This operator will be used in the regularisation term of our implicit genFRI optimisation problem. We then briefly review the method of alternating projections [31] as well as the FRI [4] framework and Cadzow denoising [21].

### A. Toeplitzification Operator

Assume that we are given an arbitrary vector $\boldsymbol{x} \in \mathbb{C}^N$, $N = 2M+1$, with entries indexed as follows:

$$\boldsymbol{x} := [x_{-M}, x_{-M+1}, \dots, x_{M-1}, x_M]^\top.$$

Then, for any $P \leq M$, we can embed $\boldsymbol{x}$ into the space $\mathbb{T}_P$ of Toeplitz matrices of $\mathbb{C}^{(N-P)\times(P+1)}$ by means of the following

[1]An efficient Python implementation of CPGD is provided on our GitHub repository [30]

[2]The alternative appellation *Toeplitzication* was used in [21].

*Toeplitzification operator:*

$$T_P : \begin{cases} \mathbb{C}^N \to \mathbb{T}_P \subset \mathbb{C}^{(N-P)\times(P+1)} \\ \boldsymbol{x} \quad \mapsto [T_P(\boldsymbol{x})]_{i,j} := x_{-M+P+i-j}, \end{cases} \quad (2)$$

where $i = 1, \ldots, N-P$, $j = 1, \ldots, P+1$. Note from (2) that the value of an entry $[T_P(\boldsymbol{x})]_{i,j}$ of the matrix $T_P(\boldsymbol{x})$ depends only on the distance $i-j$ between the row and column indices: $T_P(\boldsymbol{x})$ is therefore a *Toeplitz* matrix and the vector $\boldsymbol{x}$ is called its *generator*.

It is well-known that the Toeplitzification operator (2) can be used to implement linear convolutions. More specifically, it can be shown (see Appendix D) that the multiplication of $T_P(\boldsymbol{x})$ with a vector $\boldsymbol{u} = [u_1, \ldots, u_{P+1}]^\top \in \mathbb{C}^{P+1}$ returns the valid part[3] of the *convolution* between the two zero-padded sequences $\tilde{x} := [\ldots, 0, x_{-M}, \ldots, \boxed{x_0}, \ldots, x_M, 0, \ldots] \in \mathbb{C}^{\mathbb{Z}}$ and $\tilde{u} := [\ldots, \boxed{0}, u_1, \ldots, u_{P+1}, 0, \ldots] \in \mathbb{C}^{\mathbb{Z}}$ (following the notational convention of [1], we mark the zeroth element of a sequence $x \in \mathbb{C}^{\mathbb{Z}}$ by enclosing it in a box).

The pseudoinverse $T_P^\dagger : \mathbb{C}^{(N-P)\times(P+1)} \to \mathbb{C}^N$ of the Toeplitzification operator maps a Toeplitz matrix $\boldsymbol{H} \in \mathbb{C}^{(N-P)\times(P+1)}$ onto its generator $\boldsymbol{h} \in \mathbb{C}^N$. As shown in Appendix E, the latter is given by

$$T_P^\dagger = \Gamma^{-1} T_P^*, \quad (3)$$

where $T_P^* : \mathbb{C}^{(N-P)\times(P+1)} \to \mathbb{C}^N$ is the adjoint of the Toeplizification operator given by (see Proposition E.1)

$$T_P^* : \begin{cases} \mathbb{C}^{(N-P)\times(P+1)} \to \mathbb{C}^N \\ \boldsymbol{H} \mapsto h_j = \sum_{i=k+j-1-P} H_{ik}, \quad j = 1, \ldots, N, \end{cases} \quad (4)$$

and $\Gamma = T_P^* T_P \in \mathbb{C}^{N\times N}$ is a diagonal matrix with entries given by (see Proposition E.2):

$$\Gamma_{i,i} = \min(i, P+1, N+1-i), \quad i = 1, \ldots, N. \quad (5)$$

Observe that the composition of $T_P^*$ and $\Gamma^{-1}$ in the expression of the pseudoinverse (3) implements a *diagonal averaging:* $T_P^*$ first sums across each diagonal of the matrix $\boldsymbol{H} \in \mathbb{C}^{(N-P)\times(P+1)}$ and $\Gamma^{-1}$ then divides the sums by the number of elements on each diagonal. It is interesting to note that this operation is also leveraged in Cadzow denoising as described in [19], in order to map back the data from their high dimensional Toeplitz embedding. The formal interpretation of this diagonal averaging as the pseudoinverse of the Toeplitzification operator proposed here is nevertheless not discussed in [19], nor anywhere else we may be aware of.

## *B. FRI in a Nutshell*

The classical FRI framework, introduced in [4], aims at estimating the innovations $\{(x_k, t_k),\ k = 1, \ldots, K\} \subset \mathbb{C} \times [0, T[$, of a $T$-periodic stream of Diracs:

$$x(t) = \sum_{k' \in \mathbb{Z}} \sum_{k=1}^{K} x_k \delta(t - t_k - Tk'), \quad \forall t \in \mathbb{R}.$$

In standard FRI, the estimation procedure is divided into two stages. The locations $t_k$ are first estimated by a nonlinear method, and then arranged into a Vandermonde system whose solution yields the Dirac amplitudes [19]. The recovery of the locations $t_k$ relies on the so-called *annihilating equation*, dating from Prony's work [32], which cancels out the Fourier series coefficients of $x$ by convolving them with a particular filter, called the *annihilating filter*. The latter is defined as the finite-tap sequence $h = [\cdots, 0, \boxed{h_0}, h_1, \ldots, h_K, 0, \cdots] \in \mathbb{C}^{\mathbb{Z}}$, with $z$-transform vanishing at roots $\{u_k := e^{-j2\pi t_k/T},\ k = 1, \ldots, K\}$:

$$H(z) = \sum_{k=0}^{K} h_k z^{-k} = h_0 \prod_{k=1}^{K} (1 - u_k z^{-1}). \quad (6)$$

For such a filter, we have indeed

$$\begin{aligned}(\hat{x} * h)_m &= \sum_{k=0}^{K} h_k \hat{x}_{m-k} \\ &= \sum_{k'=1}^{K} x_{k'} \left( \sum_{k=0}^{K} h_k u_{k'}^{-k} \right) u_{k'}^m = 0,\ m \in \mathbb{Z}, \end{aligned} \quad (7)$$

where $\hat{x}_m = \sum_{k=1}^K x_k u_k^m$, $m \in \mathbb{Z}$, are the Fourier coefficients of $x$ in (1). Notice that the roots $u_k$ of the z-transform $H(z)$ in (6) of $h$ are , ignoring multiplicative constants, in one-to-one correspondence with the locations $t_k$. Recovering them amounts to estimating the coefficients $\boldsymbol{h} = [h_0, \ldots, h_K] \in \mathbb{C}^{K+1}$ of $h$ from the annihilating equation (7). If for instance we have $N = 2M+1$ consecutive Fourier coefficients of $x$, e.g. $\boldsymbol{x} = [\hat{x}_{-M}, \ldots, \hat{x}_M] \in \mathbb{C}^{2M+1}$, we can extract the $N-K$ equations from (7) corresponding to the convolution indices $m = -M+K, \ldots, M$, and use the Toeplitzification operator.[4] defined in (2) to form the following matrix equation:

$$T_K(\boldsymbol{x})\boldsymbol{h} = 0_{N-K}, \quad \|\boldsymbol{h}\| \neq 0. \quad (8)$$

Observe that any nontrivial element of the nullspace of $T_K(\boldsymbol{x})$ is a solution to (8). For $M \geq K$, it can be shown [19] that $T_K(\boldsymbol{x}) \in \mathbb{C}^{(N-K)\times(K+1)}$ has rank $K$ and therefore has a nontrivial nullspace with dimension 1. Up to a multiplicative constant, the annihilating equation (8) admits hence a unique solution. The latter can be obtained numerically by means of *total least-squares* [19], which computes the eigenvector associated to the smallest[5] eigenvalue of $T_K(\boldsymbol{x})$. In the critical case $M = K$, the matrix $T_K(\boldsymbol{x})$ is *square*, while in the *oversampling* case $M > K$ it is *rectangular and tall*. As explained in [19], oversampling makes the estimation procedure more resilient to potential noise perturbations in the Fourier coefficients. In such cases, Blu *et al.* recommend moreover to perform Cadzow denoising on the Fourier coefficients $\boldsymbol{x}$ (see Section II-C) as well as replace (8) by a more general annihilating equation:

$$T_P(\boldsymbol{x})\tilde{\boldsymbol{h}} = 0_{N-P}, \quad \|\tilde{\boldsymbol{h}}\| \neq 0, \quad (9)$$

[3]See Appendix D for a formal definition of the valid part of a convolution between zero-padded sequences.

[4]Remember the link between the Toeplitzification operator and convolution discussed in Section II-A

[5]An eigenvalue exactly equal to zero may in practice be impossible to obtain due to numerical inaccuracies.

with $K \leq P \leq M$, and $\tilde{\boldsymbol{h}} \in \mathbb{C}^{P+1}$. Again, it is possible to show that $T_P(\boldsymbol{x})$ has rank $K$, and hence a nontrivial nullspace with dimension $P+1-K$. Solutions to (9) are therefore not unique in this case, but all are equally valid for practical purposes. Moreover, the increased nullspace dimension makes Cadzow denoising more efficient at filtering the noise component. In practice, the case $P = M$ has been reported to yield the best empirical performance [19].

### C. Cadzow Denoising

For strong noise perturbations, the generalised annihilating equation (9) may fail to admit a nontrivial solution. Indeed, noisy generators $\boldsymbol{x}$ can yield full column rank matrices $T_P(\boldsymbol{x})$ with trivial nullspace. As a potential cure, Blu *et al.* propose to denoise the Fourier coefficients $\boldsymbol{x}$ prior to solving the annihilating equation. This denoising step attempts to transform $T_P(\boldsymbol{x})$ into a Toeplitz matrix with rank at most $K$, thus guaranteeing the existence of nontrivial solutions to (9). This operation is carried out by means of *Cadzow denoising* [21], an alternating projection method[6] applied heuristically to the subspace $\mathbb{T}_P$ of Toeplitz matrices and the subset $\mathcal{H}_K$ of matrices with rank at most $K$:

$$\mathcal{H}_K := \left\{\boldsymbol{M} \in \mathbb{C}^{(N-P)\times(P+1)} \mid \operatorname{rank} \boldsymbol{M} \leq K\right\}. \tag{10}$$

Using the notation introduced in Section II-A, Cadzow denoising can be seen as processing the noisy coefficients $\boldsymbol{x}$ as follows:

$$\check{\boldsymbol{x}} = T_P^\dagger \left[\Pi_{\mathbb{T}_P}\Pi_{\mathcal{H}_K}\right]^n T_P(\boldsymbol{x}), \tag{11}$$

for some suitable $n \in \mathbb{N}$, where $\Pi_{\mathbb{T}_P}$ and $\Pi_{\mathcal{H}_K}$ are the projections onto the subsets $\mathbb{T}_P$ and $\mathcal{H}_K$ of $\mathbb{C}^{(N-P)\times(P+1)}$ respectively.

Note that since $\mathcal{H}_K$ is a *non convex* set the convergence of the method of alternating projection (MAP) in (11) is not guaranteed. Nevertheless, experimental results [19], [21] suggest that Cadzow denoising almost always converges after a few iterations (typically $n \leq 20$), which could theoretically[7] be explained by the *local* convergence result Theorem G.1 discussed in Appendix G. We conclude this section by providing closed-form expressions for the projection operators $\Pi_{\mathbb{T}_P}$ and $\Pi_{\mathcal{H}_K}$, needed in (11).

*1) Projection Onto $\mathbb{T}_P$:* As shown in Appendix F, the orthogonal projection operator onto the subspace $\mathbb{T}_P \subset \mathbb{C}^{(N-P)\times(P+1)}$ of rectangular Toeplitz matrices can be written in terms of the Toeplitzification operator and its pseudoinverse as:

$$\Pi_{\mathbb{T}_P} = T_P T_P^\dagger = T_P \Gamma^{-1} T_P^*. \tag{12}$$

*2) Projection Onto $\mathcal{H}_K$:* The projection operator onto the space $\mathcal{H}_K$ of matrices with rank at most $K$ is given by the *Eckart-Young-Minsky theorem* [33]. The latter states that the projection map

$$\Pi_{\mathcal{H}_K}(\boldsymbol{X}) = \arg\min_{\boldsymbol{H}\in\mathcal{H}_K} \|\boldsymbol{X}-\boldsymbol{H}\|_F, \quad \boldsymbol{X} \in \mathbb{C}^{(N-P)\times(P+1)}, \tag{13}$$

can be computed in closed-form as:

$$\Pi_{\mathcal{H}_K}(\boldsymbol{X}) = \boldsymbol{U}\Lambda_K\boldsymbol{V}^H, \quad \boldsymbol{X} \in \mathbb{C}^{(N-P)\times(P+1)}, \tag{14}$$

where $\boldsymbol{X} = \boldsymbol{U}\Lambda\boldsymbol{V}^H$ is the singular value decomposition of $\boldsymbol{X}$, and $\Lambda_K$ is the diagonal matrix of sorted singular values truncated to the $K$ strongest ones. Note that the output of the projection map is unique as long as the $K-$th and $(K+1)-$th largest singular values are different. Fortunately, the space of matrices failing to verify this condition is very small –more precisely it is *thin*, as discussed extensively in [34, Section 2]. In practice moreover, floating-point arithmetic makes it very unlikely that the $K-$th and $(K+1)-$th largest singular values be exactly identical. Thus, the projection map $\Pi_{\mathcal{H}_K}$ can be considered single-valued for practical purposes.

## III. Generalised FRI As an Inverse Problem

### A. Generalised FRI

In Section II-B, we have described a procedure for recovering the locations $t_k$ from consecutive Fourier coefficients of $x$. The issue of computing these Fourier coefficients from a collection of arbitrary linear *measurements* $\boldsymbol{y} \in \mathbb{C}^L$ of $x$, $L \geq N$ now remains. Blu *et al.* [19] treated the simple scenario of measurements resulting from regular time sampling with ideal low-pass prefiltering. In such a case, they showed that, for a well chosen prefilter bandwidth, the Fourier coefficients could simply be obtained by applying a discrete Fourier transform to the measurements $\boldsymbol{y}$. For more general measurement types, the situation is more complex, and the Fourier coefficients $\boldsymbol{x} \in \mathbb{C}^N$ must in general be estimated by solving a linear inverse problem:

$$\boldsymbol{y} = \boldsymbol{G}\boldsymbol{x} + \boldsymbol{n}, \tag{15}$$

where the forward matrix $\boldsymbol{G} \in \mathbb{C}^{L\times N}$, $L \geq N$, is application dependent, and $\boldsymbol{n}$ is additive noise, usually assumed to be a white Gaussian random vector. In [23], Pan *et al.* proposed the *generalised FRI (genFRI)* optimisation problem to deal with (15). The latter is a non convex constrained optimisation problem whose objective is to jointly recover the Fourier coefficients $\boldsymbol{x} \in \mathbb{C}^N$ –required to minimise a quadratic data-fidelity term– and their corresponding annihilating filter coefficients $\boldsymbol{h} \in \mathbb{C}^{P+1}$. The annihilating equation linking the two unknowns is explicitly enforced as a constraint, yielding an optimisation problem of the form:

$$\min_{\substack{\boldsymbol{x}\in\mathbb{C}^N\\ \boldsymbol{h}\in\mathbb{C}^{P+1}}} \|\boldsymbol{G}\boldsymbol{x}-\boldsymbol{y}\|_2^2 \quad \text{subject to} \quad \begin{cases} T_P(\boldsymbol{x})\boldsymbol{h} = 0_{N-P},\\ \langle \boldsymbol{h}, \boldsymbol{h}_0\rangle = 1, \end{cases} \tag{16}$$

where $\boldsymbol{h}_0 \in \mathbb{C}^{P+1}$ is generated randomly according to the circularly-symmetric complex Gaussian distribution $\mathbb{C}\mathcal{N}(0, \boldsymbol{I}_{P+1})$. The normalisation constraint,[8] $\langle \boldsymbol{h}, \boldsymbol{h}_0\rangle = 1$ is used to exclude trivial solutions to the annihilating equation in (16) [13], [23].

[6] See Appendix G for a review of the method of alternating projections.

[7] As explained in Appendix G however, the assumptions of Theorem G.1 are unfortunately very difficult to verify in practice.

[8] In [23] the authors have also considered the more natural normalisation constraint $\|\boldsymbol{h}\| = 1$. They claim however that this normalisation strategy is less successful experimentally.

Pan *et al.* propose to solve (16) via an *heuristic* alternating minimisation algorithm described in Appendix H. At each iteration, the annihilating filter $\boldsymbol{h}_n$ and the Fourier series coefficients $\boldsymbol{x}_n$ are updated by solving two linear systems of size $2(N+1) \times 2(N+1)$ and $(2N-P) \times (2N-P)$ respectively. In practice, the algorithm is stopped when the data mismatch $\|\boldsymbol{G}\boldsymbol{x}_n - \boldsymbol{y}\|_2$ falls below a certain threshold $\epsilon$ (typically the noise level). Note that this heuristic iterative procedure comes without any strong or weak convergence guarantee: it is not known if the sequence $\{(\boldsymbol{h}_n, \boldsymbol{x}_n),\ n \in \mathbb{N}\} \subset \mathbb{C}^{P+1} \times \mathbb{C}^N$ converges and if its limit coincides with a critical point of (16). Moreover, the proposed stopping criterion requires a knowledge of the noise level, unknown in practice. When it is unknown, Pan *et al.* recommend performing the reconstruction for a fixed and arbitrary number of iterations (typically 50). For optimal performances, they also suggest to run the algorithm for multiple random initialisations of $\boldsymbol{h}_0 \in \mathbb{C}^{P+1}$ (typically 15). The overall reconstruction procedure can therefore be computationally intensive since each iteration has a complexity of $\mathcal{O}(8(N+1)^3 + (2N-P)^3) = \mathcal{O}(N^3)$ –the cost of solving the linear systems.

### B. Implicit Generalised FRI

The annihilating equation constraint in (16) can be thought of as regularising the genFRI problem. Indeed, minimising the quadratic term $\|\boldsymbol{G}\boldsymbol{x} - \boldsymbol{y}\|_2^2$ alone in the presence of noise would not necessarily yield Fourier coefficients $\boldsymbol{x}$ with nontrivial annihilating filter, which the annihilating constraint enforces explicitly. Unfortunately, this regularisation also complicates significantly the optimisation procedure. Indeed, it requires the introduction of an extra unknown variable with non linear dependency on the data, namely the annihilating filter $\boldsymbol{h}$. Moreover, the non linear constraint $T_P(\boldsymbol{x})\boldsymbol{h} = 0_{N-P}$ is highly non convex, and state-of-the-art algorithms, such as alternating minimisation or gradient descent [24], may suffer from getting trapped in local minima.[9] [35]. To circumvent these issues, we propose the following *implicit* formulation of the genFRI problem, in which only the Fourier coefficients are recovered:

$$\min_{\boldsymbol{x} \in \mathbb{C}^N} \|\boldsymbol{G}\boldsymbol{x} - \boldsymbol{y}\|_2^2 \quad \text{subject to} \quad \begin{cases} \operatorname{rank} T_P(\boldsymbol{x}) \leq K, \\ \|\boldsymbol{x}\|_2 \leq \rho, \end{cases} \tag{17}$$

where $K \leq P \leq M$ and $\rho \in ]0, +\infty]$.

Similarly to (16), the quadratic term $\|\boldsymbol{G}\boldsymbol{x} - \boldsymbol{y}\|_2^2$ in (17) is used to guarantee high fidelity of the recovered coefficients to the observed data. Unlike (16), (17) leverages a regularising rank constraint on $T_P(\boldsymbol{x})$ which does not explicitly involve the unknown annihilating filter. As already discussed in Section II-C in the context of Cadzow denoising, requiring $T_P(\boldsymbol{x})$ to be of rank at most $K$ is indeed a sufficient condition for the generalised annihilating equation (9) to admit nontrivial solutions. This implicit regularisation greatly simplifies the genFRI problem, since it decouples the problem of estimating the Fourier coefficients from the problem of estimating the annihilating filter. The normalisation constraint $\|\boldsymbol{x}\|_2 \leq \rho$ enforces finite energy to the recovered Fourier coefficients. As shall be seen in Section IV, it can be relaxed when the forward matrix $\boldsymbol{G}$ is injective by setting $\rho = +\infty$. Indeed, it is only used to ensure *coercivity* in underdetermined cases where the forward matrix $\boldsymbol{G}$ has a nontrivial null space. Coercivity is indeed a key assumption [36] for the convergence of the proximal gradient descent method envisioned in Section IV-A.

*Remark (On the choice of P):* Note that the rank constraint in (17) is more *selective* for values of $P$ close to $M$, hence enforcing a stronger regularisation. Indeed, it is easy to see that the maximal rank of rectangular matrices in $\mathbb{C}^{(N-P)\times(P+1)}$ ranges in[10] $[\![K+1, M+1]\!]$ when $P$ ranges in $[\![K, M]\!]$. Consequently, the subset $\mathcal{H}_K$ of matrices of rank at most $K$ becomes "smaller and smaller" relatively to the ambient space as $P$ increases towards $M$. We can hence expect (17) to perform better in practice for $P = M$. This is in contrast with the explicit generalised FRI problem (16), whose equality constraint is equally stringent for different values of $P$.

*Remark (Case $\boldsymbol{G} = \boldsymbol{I}$):* When $\boldsymbol{G} = \boldsymbol{I}$, the optimisation problem (17) becomes a simple denoising problem, which could therefore be used as an alternative to Cadzow denoising or its upgrade [21].

## IV. Optimisation Algorithm

### A. Non Convex Proximal Gradient Descent

The optimisation problem (17) can be rewritten in an unconstrained form as:

$$\min_{\boldsymbol{x} \in \mathbb{C}^N} \|\boldsymbol{G}\boldsymbol{x} - \boldsymbol{y}\|_2^2 \; + \; \iota_{\mathcal{H}_K}(T_P(\boldsymbol{x})) \; + \; \iota_{\mathbb{B}_\rho}(\boldsymbol{x}), \tag{18}$$

where $\mathcal{H}_K$ is the non convex set of matrices with rank lower than or equal to $K$ defined in (10), $\mathbb{B}_\rho := \{\boldsymbol{x} \in \mathbb{C}^N : \|\boldsymbol{x}\|_2 \leq \rho\}$ is the $\ell_2$-ball with radius $\rho > 0$, and $\iota_{\mathcal{H}_K} : \mathbb{C}^{(N-P)\times(P+1)} \to \{0, +\infty\}$, $\iota_{\mathbb{B}_\rho} : \mathbb{C}^N \to \{0, +\infty\}$ are indicator functions with domains $\mathcal{H}_K$ and $\mathbb{B}_\rho$, respectively. Observe that the unconstrained optimisation problem (18) can be written as a sum between a *convex* and *differentiable* quadratic term

$$F(\boldsymbol{x}) := \|\boldsymbol{G}\boldsymbol{x} - \boldsymbol{y}\|_2^2, \qquad \boldsymbol{x} \in \mathbb{C}^N,$$

and a *non convex* and *non differentiable* term

$$H(\boldsymbol{x}) := \iota_{\mathcal{H}_K}(T_P(\boldsymbol{x})) \; + \; \iota_{\mathbb{B}_\rho}(\boldsymbol{x}), \qquad \boldsymbol{x} \in \mathbb{C}^N.$$

It is moreover easy to see that the gradient of $F$

$$\nabla F(\boldsymbol{x}) = 2\boldsymbol{G}^H(\boldsymbol{G}\boldsymbol{x} - \boldsymbol{y}), \quad \boldsymbol{x} \in \mathbb{C}^N, \tag{19}$$

is $\beta$*-Lipschitz continuous* with Lipschitz constant given by twice the spectral norm of the matrix $\boldsymbol{G}^H\boldsymbol{G}$:

$$\beta = 2\left\|\boldsymbol{G}^H\boldsymbol{G}\right\|_2 = \sup\left\{2\left\|\boldsymbol{G}^H\boldsymbol{G}\boldsymbol{x}\right\|_2 : \boldsymbol{x} \in \mathbb{C}^N,\ \|\boldsymbol{x}\|_2 = 1\right\}. \tag{20}$$

It is hence possible to optimise (18) by means of *proximal gradient descent (PGD)* [24], an iterative method alternating between gradient and proximal steps according to the following

[9] This is notably the reason why Pan *et al.* recommend multiple random initialisations of their algorithm in [23].

[10] For $n, m \in \mathbb{Z}$, $n < m$, we denote by $[\![n, m]\!]$ the integer interval $[n, m] \cap \mathbb{Z}$.

update equation:

$$x_{k+1} \in \operatorname{prox}_{\tau H}\left(x_k - \tau \nabla F\left(x_k\right)\right), \tag{21}$$

for $k \geq 0$, $x_0 \in \mathbb{C}^N$, $\tau > 0$ and $\operatorname{prox}_{\tau H}$ defined in (22). Given a current estimate $x_k \in \mathbb{C}^N$, the update equation (21) decreases the value of the objective function (18) by selecting a *proximal point* [24] –with respect to $H$– of a target located at a distance $\tau$ from $x_k$ along the direction of steepest descent $-\nabla F(x_k)$. The operator mapping a point $x \in \mathbb{C}^N$ to its proximal points with respect to $H$ is called *proximal operator*, and is defined as [24]

$$\operatorname{prox}_{\tau H}(x) : \begin{cases} \mathbb{C}^N \to \mathcal{P}\left(\mathbb{C}^N\right), \\ x \mapsto \arg\min_{z \in \mathbb{C}^N} \frac{1}{2\tau}\|x - z\|_2^2 + H(z), \end{cases} \tag{22}$$

where $\mathcal{P}(\mathbb{C}^N)$ is the power set of $\mathbb{C}^N$, and $\tau > 0$ controls the relative importance of $H$ with respect to the squared distance to $x$. The function $H$ being non convex, the proximal operator (22) will in general return *multiple* proximal points, which can all be used interchangeably in (21). The convergence of the sequence $\{x_k\}_{k\in\mathbb{N}}$ of PGD iterates (21) towards critical points of (18) is established by the following theorem.

*Theorem 1 (Convergence of PGD for Arbitrary **G**):* Assume that $\rho \in (0, +\infty)$ in (18), and $\tau < 1/\beta$ with $\beta$ defined in (20). Then, any *limit point* $x_\star$ of the sequence $\{x_k\}_{k\in\mathbb{N}}$ generated by (21) is a *local minimum* of (18).

*Proof:* The proof of this theorem is adapted from [36, Theorem 1] and given in Appendix I.

As stated by Theorem 2 hereafter, the convergence of PGD furthermore extends to the case $\rho = +\infty$, at least for injective forward matrices $G$. Setting $\rho = +\infty$ in (17) is equivalent to dropping the energy normalisation constraint, since $\|x\|_2 \leq +\infty$ is trivially verified and hence the associated indicator function $\iota_{\mathbb{B}_\rho}$ in (18) is always null.

*Theorem 2 (Convergence of PGD for Injective **G**):* Assume that $\rho = +\infty$ in (18), $\tau < 1/\beta$ with $\beta$ defined in (20), and $G \in \mathbb{C}^{L\times N}$ in (18) is *injective*, i.e., $\ker(G) = \{0_N\}$. Then, any *limit point* $x_\star$ of the sequence $\{x_k\}_{k\in\mathbb{N}}$ generated by (21) is a *local minimum* of (18).

*Proof:* The proof of this theorem is given in Appendix I. ■

A practical implication of Theorem 2 is that, for injective forward matrices $G$, PGD applied to the following relaxed implicit genFRI problem is convergent:

$$\min_{x\in\mathbb{C}^N} \|Gx - y\|_2^2 + \iota_{\mathcal{H}_K}\left(T_P\left(x\right)\right), \tag{23}$$

where $F(x) := \|Gx - y\|_2^2$, and $H(x) := \iota_{\mathcal{H}_K}(T_P(x))$. As discussed in Section IV-B, (23) should always be favoured over (18) for injective forward matrices $G$, since solving it via PGD requires less computations at each proximal step.

### B. Cadzow PnP Gradient Descent

As seen in the previous section, PGD requires the computation of the proximal operator (22) at each iteration, which amounts to finding a minimiser to the following non convex optimisation problem:

$$\min_{z\in\mathbb{C}^N} \left\{\frac{1}{2\tau}\|x - z\|_2^2 + \iota_{\mathcal{H}_K}\left(T_P\left(z\right)\right) + \iota_{\mathbb{B}_\rho}\left(z\right)\right\}, \tag{24}$$

for some input $x \in \mathbb{C}^N$. Observe that the proximal step (24) can be seen as a generalised *projection* step, aiming to find a point as close as possible from $x$ while verifying some convex and non convex constraints specified by the indicator functions. This is formalised by Proposition 1 hereafter, which shows that solutions to (24) can be identified with those of a projection problem:

*Proposition 1:* Consider the Toeplitz matrix $W := T_P(\operatorname{diag}(\Gamma^{-1/2}))$ where $\operatorname{diag} : \mathbb{C}^{N\times N} \to \mathbb{C}^N$ is the linear operator mapping a matrix onto its diagonal and $\Gamma = T_P^* T_P \in \mathbb{C}^{N\times N}$ is the diagonal and positive definite matrix given by (5). Then, the proximal operator (22) of $H(x) := \iota_{\mathcal{H}_K}(T_P(x)) + \iota_{\mathbb{B}_\rho}(x)$, for $\rho \in ]0, +\infty]$, $\tau > 0$ and $K \leq P \leq M$ is given by

$$\operatorname{prox}_{\tau H}(x) = T_P^\dagger \Pi^W_{\mathbb{T}_P\cap\mathcal{H}_K\cap\mathbb{B}^W_\rho} T_P(x), \quad \forall x \in \mathbb{C}^N, \tag{25}$$

where

$$\mathbb{B}^W_\rho := \{X \in \mathbb{C}^{(N-P)\times(P+1)} : \|W \odot X\|_F \leq \rho\}, \tag{26}$$

and $\Pi^W_{\mathbb{T}_P\cap\mathcal{H}_K\cap\mathbb{B}_\rho}$ is the *projection operator* onto $\mathbb{T}_P \cap \mathcal{H}_K \cap \mathbb{B}^W_\rho$ with respect to the $W$*-weighted Frobenius norm:*

$$\Pi^W_{\mathbb{T}_P\cap\mathcal{H}_K\cap\mathbb{B}^W_\rho} : \begin{cases} \mathbb{C}^{(N-P)\times(P+1)} \to \mathcal{P}\left(\mathbb{C}^{(N-P)\times(P+1)}\right), \\ X \mapsto \arg\min_{Z\in\mathbb{T}_P\cap\mathcal{H}_K\cap\mathbb{B}^W_\rho} \|W \odot (X - Z)\|_F. \end{cases} \tag{27}$$

The symbol $\odot$ in (26) and (27) denotes the *Hadamard product* for matrices.

*Proof:* The proof of this proposition is given in Appendix A. ■

Equation (25) provides us with a three-step recipe for computing the proximal operator (22) associated to a vector $x \in \mathbb{C}^N$:

1) Transform the input vector $x$ into a Toeplitz matrix via the Toeplitzification operator $T_P$.
2) Project $T_P(x)$ onto $\mathbb{T}_P \cap \mathcal{H}_K \cap \mathbb{B}^W_\rho$ by solving the *weighted structured low-rank approximation (WSLRA)* problem (27).
3) Map back the projected matrix onto $\mathbb{C}^N$ using the pseudoinverse $T_P^\dagger$ of the Toeplitzification operator.

Unfortunately, simple closed-form solutions to the WSLRA problem (27) are unavailable in general, and must therefore be approximated numerically. We propose to perform such an approximation by means of the method of alternating projections (MAP) (see Appendix G):

$$\Pi^W_{\mathbb{T}_P\cap\mathcal{H}_K\cap\mathbb{B}^W_\rho} \simeq \left[\Pi^W_{\mathbb{T}_P}\Pi^W_{\mathcal{H}_K}\Pi^W_{\mathbb{B}^W_\rho}\right]^n, \tag{28}$$

where $n \in \mathbb{N}$ and where $\Pi^W_{\mathbb{T}_P}$, $\Pi^W_{\mathcal{H}_K}$ and $\Pi^W_{\mathbb{B}^W_\rho}$ are the projection operators onto $\mathbb{T}_P$, $\mathcal{H}_K$ and $\mathbb{B}^W_\rho$ with respect to the $W$-weighted Frobenius norm. As detailed in Propositions 2 and 3 hereafter, the projection operators $\Pi^W_{\mathbb{T}_P}$ and $\Pi^W_{\mathbb{B}^W_\rho}$ admit simple closed-form expressions.

*Proposition 2 (Weighted Projection onto $\mathbb{T}_P$):* Let $W$ be the Toeplitz matrix from Proposition 1. Then, the projection operator onto $\mathbb{T}_P$ with respect to the $W$-weighted Frobenius norm is given, for every $X \in \mathbb{C}^{(N-P)\times(P+1)}$, by:

$$\Pi^W_{\mathbb{T}_P}(X) = \Pi_{\mathbb{T}_P}(X), \tag{29}$$

where $\Pi_{\mathbb{T}_P}$ is the projection operator onto $\mathbb{T}_P$ with respect to the *canonical* Frobenius norm given in (12).

*Proof:* The proof of this proposition is given in Appendix J.

*Proposition 3 (Weighted Projection onto $\mathbb{B}_\rho^{\boldsymbol{W}}$):* The projection operator onto $\mathbb{B}_\rho^{\boldsymbol{W}}$ with respect to the $\boldsymbol{W}$-weighted Frobenius norm is given, for every $\boldsymbol{X} \in \mathbb{C}^{(N-P)\times(P+1)}$, by:

$$\Pi_{\mathbb{B}_\rho^{W}}^{\boldsymbol{W}}(\boldsymbol{X}) = \begin{cases} \boldsymbol{X} \text{ if } \|\boldsymbol{W} \odot \boldsymbol{X}\|_F \leq \rho, \\ \frac{\rho \boldsymbol{X}}{\|\boldsymbol{W}\odot\boldsymbol{X}\|_F} \text{ if } \|\boldsymbol{W} \odot \boldsymbol{X}\|_F > \rho. \end{cases} \quad (30)$$

*Proof:* The proof of this proposition is given in Appendix K.■

Computing the projection operator $\Pi_{\mathcal{H}_K}^{\boldsymbol{W}}$ amounts to solving a *weighted low-rank approximation (WLRA)* problem [37], [38], which is a more difficult task. Indeed, WLRA problems admit no simple closed-form solutions and must hence be solved numerically via iterative algorithms [37], [38]. This is in contrast with the unweighted low-rank approximation problem discussed in Section II-C2, whose solutions could easily be computed via a simple truncated SVD (14). A popular algorithm for solving WLRA problems is the EM-algorithm from Srebro and Jaakkola [38] , which is relatively simple, efficient, and has good empirical performances. Starting from an initial guess $\boldsymbol{X}_0 \in \mathbb{C}^{(N-P)\times(P+1)}$, the algorithm computes the projection $\Pi_{\mathcal{H}_K}^{\boldsymbol{W}}(\boldsymbol{X})$ of a matrix $\boldsymbol{X} \in \mathbb{C}^{(N-P)\times(P+1)}$ via the following iterative scheme:

$$\boldsymbol{X}_j = \Pi_{\mathcal{H}_K}\left(\boldsymbol{W} \odot \boldsymbol{X} + (\boldsymbol{I} - \boldsymbol{W}) \odot \boldsymbol{X}_{j-1}\right), \qquad j \geq 1, \quad (31)$$

where $\Pi_{\mathcal{H}_K}$ is the projection operator (14) with respect to the unweighted Frobenius norm. Numerical experiments carried out by Srebro and Jaakkola reveal that, for weighting matrices $\boldsymbol{W}$ with non zero entries, choosing $\boldsymbol{X}_0 = \boldsymbol{X}$ promotes convergence of the iterations (31) towards a global –or at least deep local– minimum of the WLRA problem defining $\Pi_{\mathcal{H}_K}^{\boldsymbol{W}}(\boldsymbol{X})$ [38]. We therefore propose to approximate $\Pi_{\mathcal{H}_K}^{\boldsymbol{W}}$ via (31) with $\boldsymbol{X}_0 = \boldsymbol{X}$ and $j = 1$, yielding

$$\Pi_{\mathcal{H}_K}^{\boldsymbol{W}}(\boldsymbol{X}) \simeq \Pi_{\mathcal{H}_K}(\boldsymbol{W} \odot \boldsymbol{X} + (\boldsymbol{I} - \boldsymbol{W}) \odot \boldsymbol{X}) = \Pi_{\mathcal{H}_K}(\boldsymbol{X}).$$

The decision of stopping the algorithm after a single iteration is entirely motivated by computational and speed considerations. Indeed, the inexact MAP proximal step (28) –computed at each iteration of the PGD algorithm– requires $n$ successive computations of the operator $\Pi_{\mathcal{H}_K}^{\boldsymbol{W}}$ which must hence be relatively fast to compute. Since each iteration of the scheme (31) involves the computation of an (expensive) SVD, performing more than one iteration would be impractical in our context.

To summarize, the MAP (28) can be approximated in practice as:

$$\Pi_{\mathbb{T}_P\cap\mathcal{H}_K\cap\mathbb{B}_\rho^W}^{\boldsymbol{W}} \simeq \left[\Pi_{\mathbb{T}_P}^{\boldsymbol{W}}\Pi_{\mathcal{H}_K}^{\boldsymbol{W}}\Pi_{\mathbb{B}_\rho^W}^{\boldsymbol{W}}\right]^n \simeq \left[\Pi_{\mathbb{T}_P}\Pi_{\mathcal{H}_K}\Pi_{\mathbb{B}_\rho^W}^{\boldsymbol{W}}\right]^n, \quad (32)$$

where $\Pi_{\mathbb{B}_\rho^W}^{\boldsymbol{W}}$ is given by (30). Observe that when $\rho = +\infty$ (which is possible for injective matrices $\boldsymbol{G}$, see Theorem 2) we have $\Pi_{\mathbb{B}_\rho^W}^{\boldsymbol{W}} = \mathrm{Id}$ and the right-hand side of (32) simplifies to $[\Pi_{\mathbb{T}_P}\Pi_{\mathcal{H}_K}]^n$. Plugging (32) into (25) finally yields the following approximate proximal step:

$$\operatorname{prox}_{\tau H}(\boldsymbol{x}) \simeq T_P^\dagger\left[\Pi_{\mathbb{T}_P}\Pi_{\mathcal{H}_K}\Pi_{\mathbb{B}_\rho^W}^{\boldsymbol{W}}\right]^n T_P(\boldsymbol{x}), \quad \forall \boldsymbol{x} \in \mathbb{C}^N, \quad (33)$$

for some $n \geq 0$. The PGD algorithm with approximate proximal step (33) is provided in Algorithm 1. Observe that when $\rho = +\infty$, (33) reduces to Cadzow denoising (11). The effect

**Algorithm 1:** Cadzow PnP Gradient Descent (CPGD).

**Require:** $\boldsymbol{y}$, $\boldsymbol{G}$, $T_P$, $\boldsymbol{x}_0$, $K \leq P$, $\tau$ as in (36), $n \in \mathbb{N}$, $\rho > 0$
 k:=0
 **repeat**
  $z_{k+1} := \boldsymbol{x}_k - 2\tau \boldsymbol{G}^H(\boldsymbol{G}\boldsymbol{x}_k - \boldsymbol{y})$
  **if** $\rho = +\infty$ **then**
   $\boldsymbol{x}_{k+1} := T_P^\dagger[\Pi_{\mathbb{T}_P}\Pi_{\mathcal{H}_K}]^n T_P(z_{k+1})$
  **else**
   $\boldsymbol{x}_{k+1} := T_P^\dagger[\Pi_{\mathbb{T}_P}\Pi_{\mathcal{H}_K}\Pi_{\mathbb{B}_\rho^W}^{\boldsymbol{W}}]^n T_P(z_{k+1})$
  **end if**
  $k \leftarrow k + 1$
 **until**a stopping criterion is satisfied
 **return** $\boldsymbol{x}^{(k)}$

of heuristic (32) is hence to replace the proximal step in the PGD iterations by a generic *denoising* step. Such an approach is reminiscent of the *plug-and-play (PnP)* framework [28], [29] from image processing, which leverages generic denoisers to approximate complex proximal operators [27]. For this reason, we baptise our algorithm *Cadzow PnP Gradient Descent (CPGD)*. In the next section, we study the convergence of Algorithm 1.

*Remark* Note that, since $\mathcal{H}_K$ is non convex and $\Pi_{\mathcal{H}_K}^{\boldsymbol{W}}$ is approximated by $\Pi_{\mathcal{H}_K}$, the MAP approximation (32) is not guaranteed to converge towards the actual projection map (27) as $n$ grows to infinity. Empirical evidence suggests however that (32) is a good enough proxy for practical purposes. Indeed, the main role of the projection operator (27) is to regularise the PGD iterates by enforcing the structure $T_P(\boldsymbol{x}) \in \mathbb{T}_P \cap \mathcal{H}_K \cap \mathbb{B}_\rho^{\boldsymbol{W}}$ –as per the constraints of the implicit genFRI problem (17). Most often, such a behaviour is also achieved by the inexact proximal step (33). For the specific case $\rho = +\infty$, it is indeed possible to apply Theorem G.1 and show the local convergence of the MAP (32) towards a point in $\mathbb{T}_P \cap \mathcal{H}_K$ (see Corollary G.1 for a proof). This explains notably why the upgraded Cadzow algorithm from Condat and Hirabayashi [21] –which attempts to solve for the WSLRA problem directly via a heuristic primal-dual splitting method– achieves marginal accuracy gain in comparison to the more naive MAP scheme leveraged by standard Cadzow denoising.

### *C. Local Fixed-Point Convergence of CPGD*

In Section IV-A, we established Theorems 1 and 2 which show the convergence of PGD towards critical points of (18). However, such results required the computation of exact proximal steps (22) in the PGD iterations, and do not apply to CPGD which leverages the inexact proximal step (33). Convergence of PGD in non convex setups with inexact proximal steps was studied in [26], [39]. The results established in both papers require the proximal step approximation errors incurred at each iteration to be decreasing and summable, which may not necessarily be the case for the MAP approximation (32). It is nevertheless possible to demonstrate that the iterations of CPGD are locally *contractive*, and therefore locally convergent towards a fixed point using the Banach contraction principle. Such a result is stated in Theorem 3 hereafter.

*Theorem 3 (CPGD is a Local Contraction):* Let $\mathcal{R}_K \subset \mathbb{C}^{(N-P)\times(P+1)}$ be the set of matrices of rank exactly $K \leq P \leq \lfloor N/2 \rfloor$, and $U_{\tau,n} : \mathbb{C}^N \to \mathbb{C}^N$ the update CPGD map

$$U_{\tau,n}(\boldsymbol{x}) := H_n\left(\boldsymbol{x} - \tau \nabla F(\boldsymbol{x})\right), \qquad \boldsymbol{x} \in \mathbb{C}^N, \tag{34}$$

with $H_n(\boldsymbol{x}) := T_P^\dagger[\Pi_{\mathbb{T}_P}\Pi_{\mathcal{H}_K}\Pi^{\boldsymbol{W}}_{\mathbb{B}^{\boldsymbol{W}}_\rho}]^n T_P(\boldsymbol{x})$. Let $\boldsymbol{G} \in \mathbb{C}^{L\times N}$ be injective, and define

$$\beta := 2\lambda_{\max}\left(\boldsymbol{G}^H\boldsymbol{G}\right), \quad \alpha := 2\lambda_{\min}\left(\boldsymbol{G}^H\boldsymbol{G}\right),$$

where $\lambda_{\min}(\boldsymbol{M})$ and $\lambda_{\max}(\boldsymbol{M})$ denote the minimum and maximum eigenvalues of a matrix $\boldsymbol{M}$ respectively.

Then, $U_{\tau,n}$ is locally well-defined (single-valued) and *Lipschitz continuous*

$$\left\|U_{\tau,n}(\boldsymbol{x}) - U_{\tau,n}(\boldsymbol{z})\right\|_2 \leq \widetilde{L}_\tau \|\boldsymbol{x} - \boldsymbol{z}\|_2,$$

for all $\boldsymbol{x}, \boldsymbol{z} \in \mathbb{C}^N$ such that $T_P(\boldsymbol{x})$ and $T_P(\boldsymbol{z})$ are in some neighbourhood[11] of some matrix $\boldsymbol{R} \in \mathcal{R}_K$. The *Lipschitz constant* $\widetilde{L}_\tau$ is moreover independent of the local neighbourhood and given by

$$\widetilde{L}_\tau = \sqrt{P+1}L_\tau, \tag{35}$$

where $L_\tau = \max\{|1-\tau\alpha|, |1-\tau\beta|\}$. Moreover, $U_{\tau,n}$ is *contractive*, i.e., $0 \leq \widetilde{L}_\tau < 1$, for

$$\frac{1}{\beta}\left(1 - \frac{1}{\sqrt{P+1}}\right) < \tau < \frac{1}{\beta}\left(1 + \frac{1}{\sqrt{P+1}}\right). \tag{36}$$

*Proof:* The proof of this theorem is given in Appendix B. ■

The following corollary shows the local convergence of CPGD towards a fixed-point of the update map (34):

*Corollary 1 (CPGD Converges Locally):* With the same notations as in Theorem 3, assume that all CPGD iterates $\{\boldsymbol{x}_k\}_{k\in\mathbb{N}}$ are such that

$$\{T_P(\boldsymbol{x}_{k+1}), T_P(\boldsymbol{x}_k)\} \subset \mathcal{U}_k, \quad \forall k \in \mathbb{N}, \tag{37}$$

for some neighbourhoods $\{\mathcal{U}_k\}_{k\in\mathbb{N}}$ of some matrices $\{\boldsymbol{R}_k\}_{k\in\mathbb{N}} \subset \mathcal{R}_K$. Assume further that $\tau$ satisfies (36). Then, $\boldsymbol{x}_k \overset{k\to\infty}{\to} \boldsymbol{x}_\star$ where $\boldsymbol{x}_\star \in \mathbb{C}^N$ is a *fixed-point* of $U_{\tau,n}$, i.e., $U_{\tau,n}(\boldsymbol{x}_\star) = \boldsymbol{x}_\star$. Moreover, we have

$$\|\boldsymbol{x}_\star - \boldsymbol{x}_k\|_2 \leq \frac{\widetilde{L}_\tau^k}{1-\widetilde{L}_\tau}\|\boldsymbol{x}_1 - \boldsymbol{x}_0\|_2, \quad \forall k \geq 1, \tag{38}$$

where $\widetilde{L}_\tau < 1$ is the Lipschitz constant (35) of $U_{\tau,n}$.

*Proof:* The proof of this corollary is given in Appendix C. ■

*Remark (Fixed Points vs. Critical Points):* Note that Corollary 1 is a much weaker result than Theorems 1 and 2. Indeed, Corollary 1 only shows the local convergence of CPGD towards fixed points of $U_{\tau,n}$, which may not necessarily be critical points of the optimisation problem (18). Theorems 1 and 2 on the other hand, show the global convergence of PGD with exact proximal step towards critical points of (18). This is however the price to pay for computing the proximal step (24) efficiently in practice.

*Remark (Geometric Interpretation of Condition (37)):* Roughly speaking, Corollary 1 guarantees the convergence of CPGD towards a fixed point of the update map (34), provided that the forward matrix $\boldsymbol{G}$ is injective, and that any two consecutive lifted estimates $T_P(\boldsymbol{x}_k)$, $T_P(\boldsymbol{x}_{k+1})$, are in a common neighbourhood $\mathcal{U}_k$ of some matrix $\boldsymbol{R}_k \in \mathcal{R}_K$. Note that this is much less stringent than requiring the entire lifted path $\{T_P(\boldsymbol{x}_k)\}_{k\in\mathbb{N}}$ to belong to some neighbourhood $\mathcal{U}$ of some fixed matrix $\boldsymbol{R} \in \mathcal{R}_K$. Indeed, condition (37) allows the lifted estimates to travel from one neighbourhood of the manifold $\mathcal{R}_K$ to another, provided that every visited neighbourhood contains at least two consecutive lifted estimates (see Fig. 1 for an illustration). This condition, although difficult to verify in practice, seems however likely to hold for $\rho = +\infty$, small enough step sizes, large enough $n$ and $\boldsymbol{x}_0 = 0_N$. Indeed, in such a case, we have:

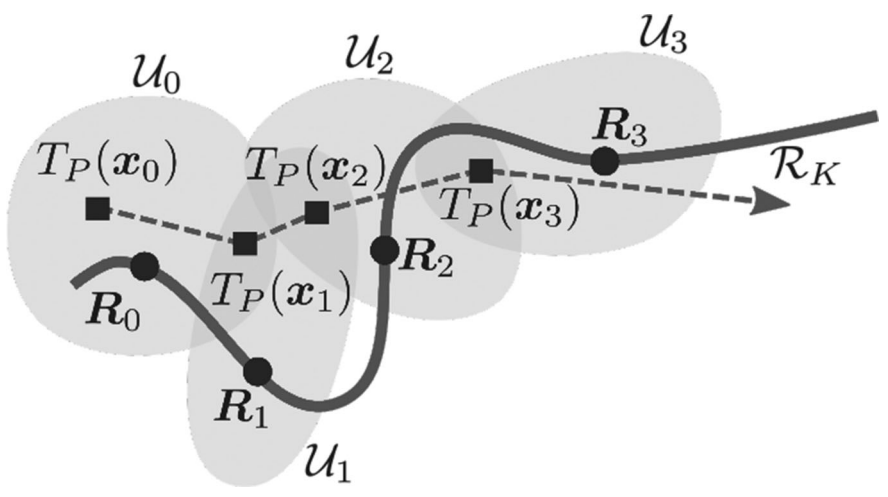

Fig. 1. Illustration of condition (37) in Corollary 1.

- $T_P(\boldsymbol{x}_0) \in \mathcal{H}_K$ is in some neighbourhood of $\mathcal{R}_K$ since $\mathcal{R}_K$ is dense in $\mathcal{H}_K$.
- For $n$ large enough, $T_P(\boldsymbol{x}_k)$ is very likely to be in some neighbourhood of $\mathcal{R}_K$, since the denoising step in the update map (34) makes $T_P(\boldsymbol{x}_k)$ close to be in the intersection $\mathcal{H}_K \cap \mathbb{T}_P$ (see Corollary G.1).
- For a small enough step size $\tau$, $T_P(\boldsymbol{x}_k)$ and $T_P(\boldsymbol{x}_{k+1})$ are likely to belong to the same neighbourhood of $\mathcal{R}_K$.

## V. Experiments and Results

In this section we validate the CPGD method numerically, considering as a testbed the scenario of irregular time sampling from [23, Section IV.A]. We assess both the reconstruction accuracy and the computational complexity of the method, and compare it to the state-of-the-art.

*Remark (Reproducibility):* Special care has been taken into making the experiments and simulations of this section fully reproducible. To reproduce the results, the reader is referred to the routines provided in our GitHub repository [30].

### A. Reconstruction Accuracy

We define a 1-periodic stream of $K = 9$ Diracs (see Fig. 2):

$$x(t) = \sum_{m\in\mathbb{Z}}\sum_{k=1}^{K} x_k\delta(t - t_k - m), \quad \forall t \in \mathbb{R}, \tag{39}$$

where the amplitudes $x_k \in \mathbb{R}_+$ and locations $t_k \in [0, 1)$ are random, with log-normal and uniform[12] distributions respectively.

[11] If $X$ is a topological space and $p$ is a point in $X$, a *neighbourhood* of $p$ is a subset $V$ of $X$ that includes an open set $U$ containing $p$, $p \in U \subseteq V$.

[12] To avoid degenerate cases, the Diracs are required to have a minimum separation distance of 1% of the total period.

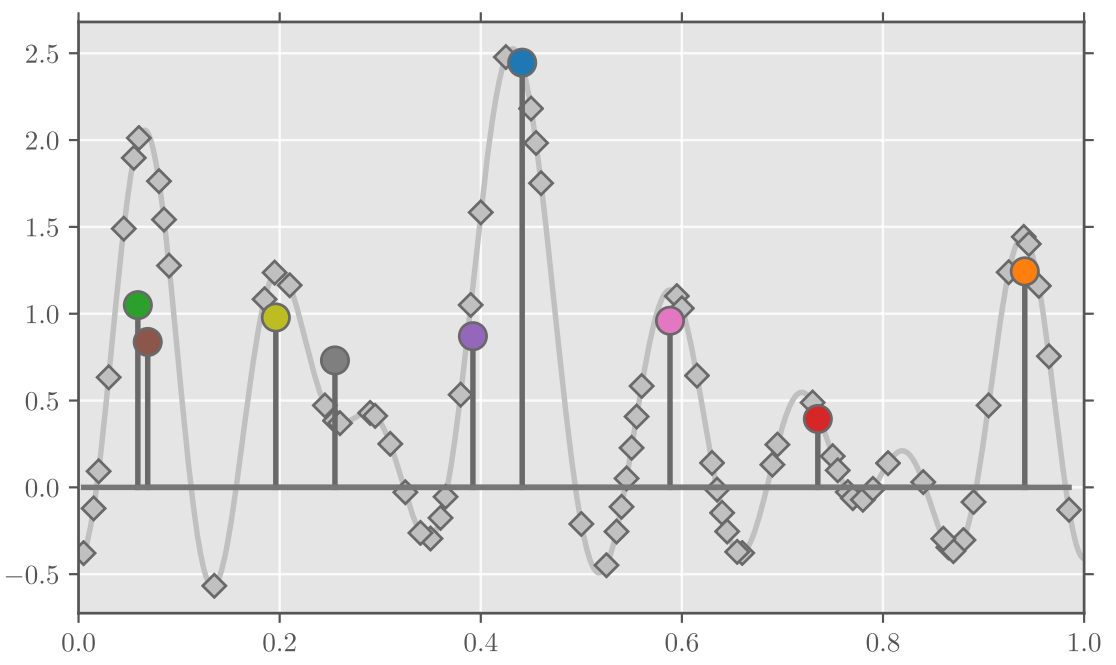

Fig. 2. Dirac stream with $K = 9$ sources (dark grey, round coloured heads) and noiseless irregular time samples (light grey diamond heads) of the low-pass filtered Dirac stream (light grey plain line), for a bandwidth $2M + 1 = 19$.

We then generate $L = 73$ noisy samples as

$$y_l = \sum_{m=-M}^{M} \hat{x}_m e^{j2\pi m\theta_l} + \epsilon_l, \quad l = 1, \ldots, L, \tag{40}$$

where $\hat{x}_m = \sum_{k=1}^{K} x_k \exp(-j2\pi m t_k)$, $m = -M, \ldots, M$, are the Fourier coefficients of the Dirac stream $x$, $\theta_l \in [0, 1)$ are chosen uniformly[13] at random, and $\epsilon_l \in \mathbb{R}$ are independent realisations of a Gaussian random variable[14] $\mathcal{N}(0, \sigma^2)$, for $l \in [\![1, L]\!]$. As explained in [23, Section IV.A], the measurements $y_l$ correspond to noisy samples of the low-pass filtered Dirac stream $x$ at irregular times $\theta_l$ (see Fig. 2), where the low-pass filter is chosen as an *ideal* low-pass filter with bandwidth $2M + 1$:

$$y_l = \sum_{k=1}^{K} x_k \varphi_M(\theta_l - t_k) \; + \; \epsilon_l, \quad l = 1, \ldots, L,$$

where

$$\varphi_M(t) := \frac{\sin\left((2M+1)\pi t\right)}{(2M+1)\sin\left(\pi t\right)}, \quad \forall t \in [0, 1),$$

is the 1-periodic sinc function or Dirichlet kernel. Using the formalism of Section III, we can rewrite (40) in vector notation as

$$\boldsymbol{y} = \boldsymbol{G}\boldsymbol{x} + \boldsymbol{\epsilon}, \tag{41}$$

where $\boldsymbol{y} = [y_1, \ldots, y_L] \in \mathbb{C}^L$, $\boldsymbol{x} = [\hat{x}_{-M}, \ldots, \hat{x}_M] \in \mathbb{C}^{N=2M+1}$, $\boldsymbol{\epsilon} = [\epsilon_1, \ldots, \epsilon_L] \in \mathbb{C}^L$, and $\boldsymbol{G} \in \mathbb{C}^{L\times N}$ is given by

$$\boldsymbol{G} = \begin{bmatrix} e^{-j2\pi M\theta_1} & \cdots 1 \cdots & e^{j2\pi M\theta_1} \\ e^{-j2\pi M\theta_2} & \cdots 1 \cdots & e^{j2\pi M\theta_2} \\ \vdots & \cdots \vdots \cdots & \vdots \\ e^{-j2\pi M\theta_{L-1}} & \cdots 1 \cdots & e^{j2\pi M\theta_{L-1}} \\ e^{-j2\pi M\theta_L} & \cdots 1 \cdots & e^{j2\pi M\theta_L} \end{bmatrix}.$$

Note that from the periodicity of complex exponentials, it is possible to flip the columns of $\boldsymbol{G}$ so as to rewrite it as a Vandermonde matrix [23]. This shows that $\boldsymbol{G}$ is *injective* when $2M + 1 \leq L$ and the irregular time samples are all distinct. From the samples $\boldsymbol{y}$ and the data model (41), we consider recovering the Fourier coefficients $\boldsymbol{x} \in \mathbb{C}^N$ by means of three algorithms:

- The CPGD algorithm 1 with $\rho = +\infty$ when $2M + 1 \leq L$ (since $\boldsymbol{G}$ is then injective) and $\rho = \|\boldsymbol{y}\|_2$ when $2M + 1 > L$ (since $\boldsymbol{G}$ is then non injective). The step size is set as $\tau = 1/\beta$ –where $\beta$ is set as described in Theorem 3– which satisfies condition (36) for the local fixed-point convergence of the algorithm (at least for injective forward matrices $\boldsymbol{G}$).
- The state-of-the-art algorithm of Pan *et al.* [23], described in Section III-A and referred to hereafter as GenFRI. For smooth integration, the Python 3 implementation of GenFRI provided by Pan *et al.* on their official Github repository [40] was included in our own algorithmic interface. Since the noise level is assumed to be unknown, we set –as recommended in [23, Section III-C]– the number of inner iterations and random initialisations to their default values, 50 and 15 respectively. Note that since GenFRI is only defined for injective matrices $\boldsymbol{G}$ (see discussion in Appendix H), we could not apply it to experimental setups with $2M + 1 > L$.
- The baseline method, referred to hereafter as LS-Cadzow, which consists in applying Cadzow denoising to the least-squares estimate of the Fourier coefficients

$$\begin{cases} \boldsymbol{x}_{\text{LS}} = \underset{\boldsymbol{x}\in\mathbb{C}^N}{\operatorname{argmin}} \|\boldsymbol{G}\boldsymbol{x} - \boldsymbol{y}\|_2^2 , \\ \boldsymbol{x}_{\text{LS-Cadzow}} = T_P^{\dagger}\left[\Pi_{\mathbb{T}_P}\Pi_{\mathcal{H}_K}\right]^n T_P\left(\boldsymbol{x}_{\text{LS}}\right). \end{cases} \tag{42}$$

We solve the least-squares optimisation problem in (42) by means of the `lstsq` function in `scipy` [41], with cut-off ratio $\texttt{cond} = 10^{-4}$.

For CPGD, we fix the maximum number of iterations[15] to 500 and consider that convergence is reached if the iterate norm is changed by less than 0.01% between two iterations. For Cadzow denoising, we fix the number of iterations to 10 for both LS-Cadzow and CPGD. For all three algorithms finally, we choose $P = M$.

The reconstruction accuracy is assessed by matching the true Dirac locations $t_k$ to the recovered ones, denoted by $\omega_k$, for $k$ between 1 and $K$. To do so, we proceed as explained in Section II-B and infer the Dirac locations $\omega_k$ from the z-transform roots of the annihilating filter associated to the Fourier coefficients estimated by each method.[16] Then, we solve the following matching problem by means of the *Hungarian algorithm*[17] [42]

$$\min_{j_1,\ldots,j_K\in\{1,\ldots,K\}} \left\{ \frac{1}{K}\sum_{k=1}^{K} d(t_k, \omega_{j_k}) \right\}, \tag{43}$$

where $d(t,\omega) = \min\{|t-\omega|, 1-|t-\omega|\}$, $\forall t,\omega \in [0,1)$, is the canonical distance on the periodised interval $[0, 1)$. Finally, we report the average positioning error, corresponding to the value of the cost function $\sum_{k=1}^{K} d(t_k, \omega_{i_k})/K$ for the indices

[13] To avoid degenerate cases, the sampling locations are required to have a minimal separation distance of 0.5% of the total period.

[14] The *noise level* is defined as the standard deviation $\sigma$ of the Gaussian distribution.

[15] In practice this upper bound is never reached: CPGD almost always converges in less than 150 iterations.

[16] See [21, Fig. 2] for additional details on the procedure used to recover the Dirac locations from the annihilating filter coefficients.

[17] The Hungarian algorithm is implemented in the `linear_sum_assignment` function from `scipy` [41].

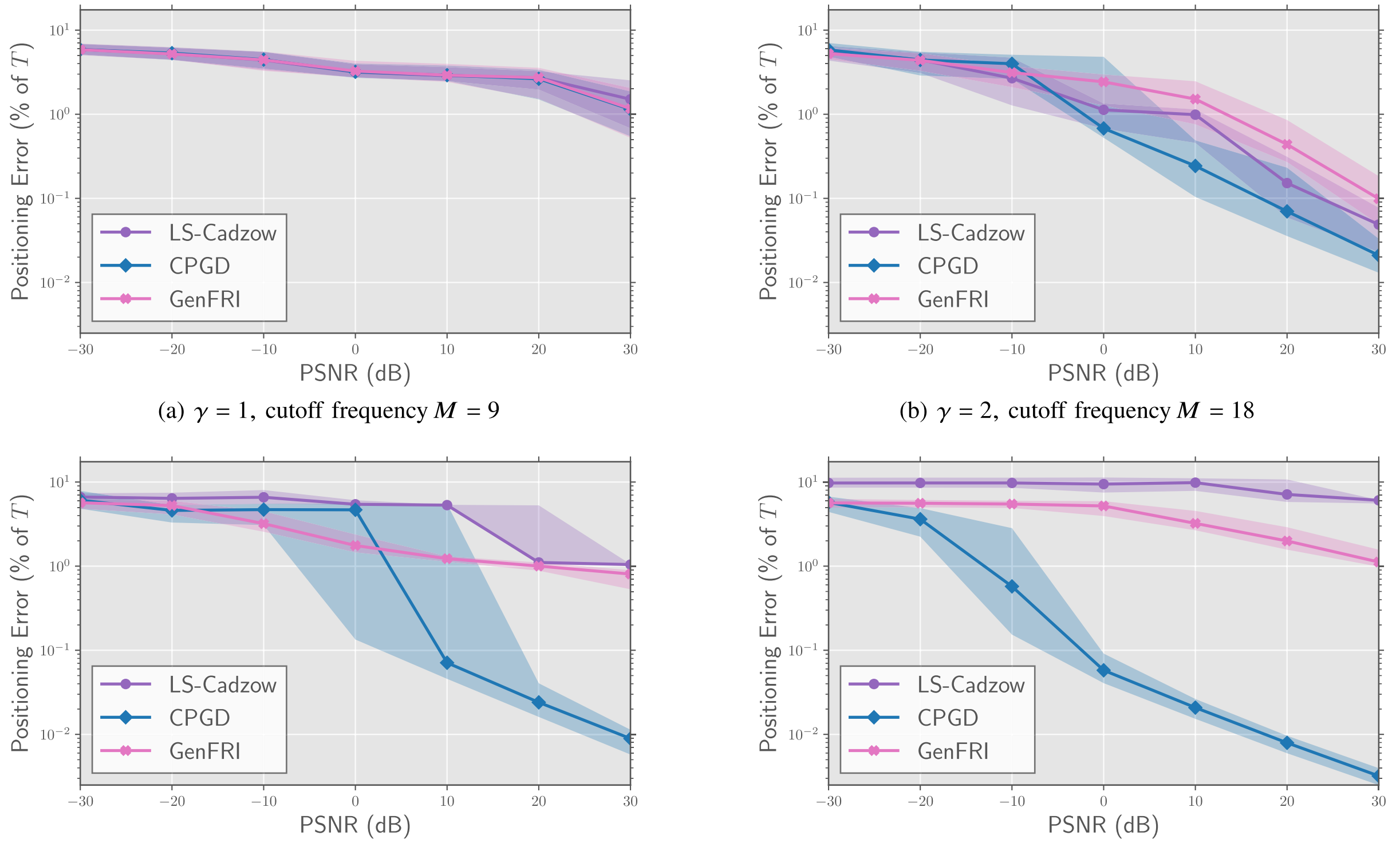

(a) $\gamma = 1$, cutoff frequency $M = 9$

(b) $\gamma = 2$, cutoff frequency $M = 18$

(c) $\gamma = 3$, cutoff frequency $M = 27$

(d) $\gamma = 4$, cutoff frequency $M = 36$

Fig. 3. Positioning error (43) (in percent of period and log-scale) for LS-Cadzow, CPGD and GenFRI, various oversampling parameters $\gamma \in \{1, 2, 3, 4\}$ and a PSNR in $\{-30, -20, -10, 0, 10, 20, 30\}$ dB. For each case, plain lines and shaded areas represent respectively the median and inter-quartile region of the positioning error's empirical distribution obtained from 192 independent noise realisations. These results can be reproduced using the Python script `reproduce_simulation_results.py` located in the directory `./benchmarking/` of our GitHub repository [30].

TABLE I
CONDITIONING NUMBER OF $G^H G$ FOR VARIOUS VALUES OF THE OVERSAMPLING PARAMETER $\gamma$

| | $\gamma = 1$ | $\gamma = 2$ | $\gamma = 3$ | $\gamma = 4$ | $\gamma = 5$ |
|---|---|---|---|---|---|
| $\kappa(G^H G)$ | 5.9 | $1.7\times10^2$ | $1.9\times10^5$ | $1.6\times10^{16}$ | $+\infty$ |

$\{i_1, \dots, i_K\}$ solutions to the matching problem (43). This metric is computed for 192 noise realisations, various cutoff frequencies $M = \gamma K$ with the oversampling factor $\gamma \in \{1, 2, 3, 4, 5\}$ (see Figs. L.1a, L.1b, L.1c, L.1d and L.1e respectively) and various noise levels

$$\sigma = \max_{k=1,\dots,K} |x_k| \times \exp\left(-\frac{\text{PSNR}}{10}\right),$$

where the peak signal to noise ratio PSNR ranges from $-30$ dB to 30 dB. The conditioning numbers of the matrix $G^H G$ for the different values of the oversampling parameter $\gamma$ are provided in Table I. The results of the experiments are displayed on Figs. 3, 4, L.2 and L.3. In Figs. 3 and 4 we plot, for different oversampling factors and PSNR, the median and inter-quartile region of the empirical distribution of the average positioning error of the different methods. In Figs. L.2 and L.3, we plot, for each source, different oversampling factors and PSNR, the median of the empirical distribution of the source location as estimated by the three methods against the true source location. The conclusions that can be drawn from the results are the following:

- Figs. 3 a, L.2a, L.2b and L.2c reveal that without oversampling in the Fourier domain (i.e., $\gamma = 1$ and a minimal cutoff frequency of $M = K = 9$), the three algorithms CPGD, GenFRI and LS-Cadzow perform similarly throughout the entire PSNR range. The average positioning error –almost indistinguishable for the three algorithms– goes from approximately 10% of the period for PSNRs of $-30$ dB to 1% of the period for PSNRs of 30 dB.

- Figs. 3 b, L.2d, L.2e and L.2f reveal reveal that with an oversampling of $\gamma = 2$ –yielding a cutoff frequency of $M = 18$, the three algorithms CPGD, GenFRI and LS-Cadzow start behaving differently for PSNRs larger than 0 dB: CPGD has the lowest positioning error, followed by LS-Cadzow and finally GenFRI. The inter-quartile regions of the positioning error's distributions are however overlapping, which means that the differences in performance are not statistically significant. For PSNRs larger than 0 dB, all algorithms have a lower positioning error than in the case $\gamma = 1$. For high PSNRs, this improvement can be as high as one and a half order of magnitude.
- Figs. 3 c, L.2g, L.2h and L.2i reveal that with an oversampling of $\gamma = 3$ –yielding a cutoff frequency of $M = 27$, the three algorithms CPGD, GenFRI and LS-Cadzow again behave differently for PSNRs larger than 0 dB: CPGD has the lowest positioning error, followed by GenFRI and finally LS-Cadzow. For high PSNRs, the differences in performance among the three algorithms become statistically significant: the inter-quartile regions of the positioning error's empirical distribution do not overlap anymore. CPGD moreover reaches a positioning error as low as 0.01% of the period, which is up to two orders of magnitude smaller than the minimal positioning error of GenFRI or LS-Cadzow in this scenario. For PSNRs greater than 10 dB, CPGD improves its positioning error with respect to the case $\gamma = 2$ by a bit less than half an order of magnitude. This is not the case for GenFRI and LS-Cadzow which both underperform with respect to the case $\gamma = 2$ –and even with respect to the case $\gamma = 1$ for LS-Cadzow. This can be explained by the large conditioning number of the matrix $\boldsymbol{G}^H\boldsymbol{G}$ in this case (see Table I), affecting the numerical stability of both algorithms (see remark in Appendix H of the suppplementary material).
- Figs. 3 a, L.2j, L.2k and L.2l reveal that with an oversampling of $\gamma = 4$ –yielding a critical bandwidth of $2M + 1 = 73$ equal to the number of measurements $L$, CPGD is superior to GenFRI which is itself superior to the baseline method LS-Cadzow in nearly all cases, with the exception of very low PSNRs ($\sim -30$ dB), where the three methods have comparable reconstruction accuracy. For PSNRs larger than $-20$ dB, the differences in performance are statistically significant. CPGD is more accurate than GenFRI and LS-Cadzow by a few orders of magnitude (from 1 to 3 orders of magnitude for PSNRs larger than $-10$ dB), reaching a minimal positioning error as low as 0.005% of the period. For PSNRs greater than $-20$ dB, CPGD improves its positioning error with respect to all previous cases $\gamma \in \{1, 2, 3\}$. Again, this is not the case for GenFRI and LS-Cadzow which both perform as good as or worse than the case $\gamma = 1$. This can be explained by the (very) large conditioning number of the matrix $\boldsymbol{G}^H\boldsymbol{G}$ in this case (see Table I), which severely affects the numerical stability of both algorithms.
- Figs. 4, L.3a, L.3b and L.3c reveal that with an oversampling of $\gamma = 5$ –yielding a critical bandwidth of $2M + 1 = 91$ greater than the number of measurements $L$, CPGD is superior to the baseline method LS-Cadzow[18] for PSNRs smaller or equal to 0 dB. The differences in reconstruction accuracy are moreover statistically significant: the positioning error of CPGD is between half an order and one order of magnitude smaller than the one of LS-Cadzow. For PSNRs between 10 and 20 dB, the two methods have however comparable performances. For PSNRs of 30 dB finally, LS-Cadzow outperforms CPGD by two orders of magnitudes. This is because when $\boldsymbol{G}$ is fat, the LS estimate $\boldsymbol{x}_{LS}$ in the LS-Cadzow algorithm (42) perfectly matches the data –i.e., $\boldsymbol{G}\boldsymbol{x}_{LS} = \boldsymbol{y}$– which is desirable for very high PSNRs. It is interesting to observe that, in contrast with the previous reconstruction accuracy profiles, CPGD's positioning error is in this case non monotonically decreasing with respect to the PSNR. This is a surprising behaviour, for which we do not have any satisfying explanation yet. Finally, it shall be noted that despite the convergence of CPGD being not proven in this case (Corollary 1 is for injective matrices $\boldsymbol{G}$ only), we did not experience any convergence issues during our numerical experiments –CPGD always converged in less than 150 iterations. Moreover, the value of the parameter $\rho > 0$ seemed to have a limited effect on the reconstruction accuracy of CPGD.

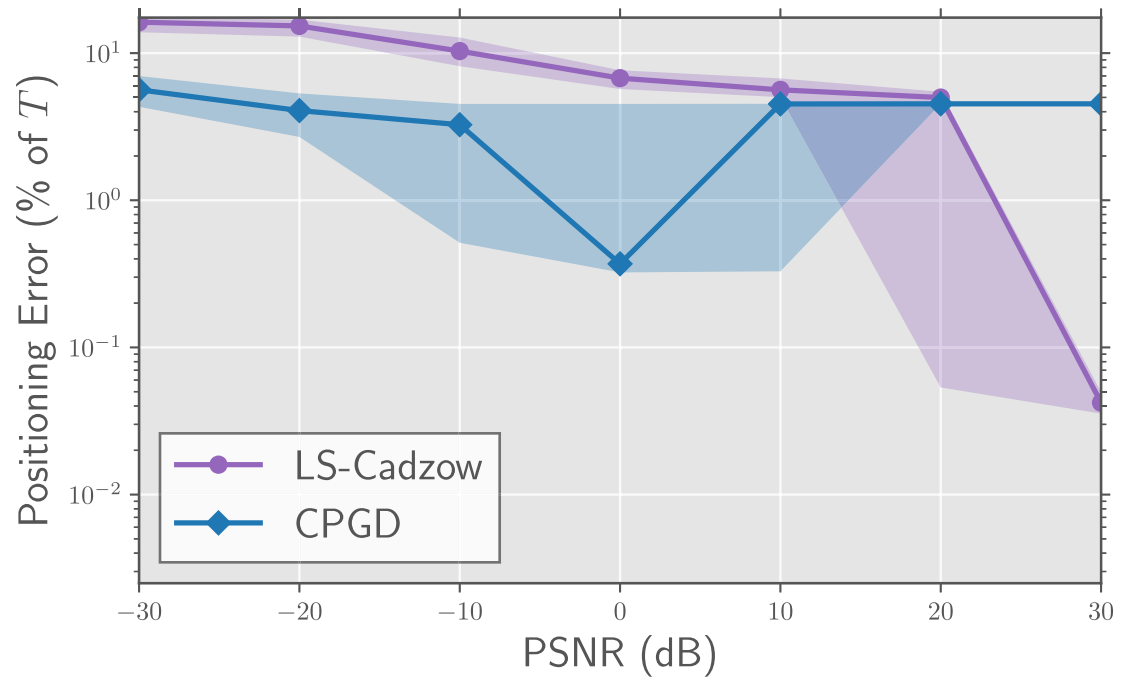

Fig. 4. Positioning error (43) (in percent of period and log-scale) for LS-Cadzow and CPGD for $\gamma = 5$ and a PSNR in $\{-30, -20, -10, 0, 10, 20, 30\}$ dB. For each case, plain lines and shaded areas represent respectively the median and inter-quartile region of the positioning error's empirical distribution obtained from 192 independent noise realisations.

In conclusion, in these simulations CPGD is better at leveraging oversampling in the Fourier domain to improve the reconstruction accuracy by several orders of magnitude with respect to the non oversampled case. In particular, CPGD performs best when the bandwidth of the low-pass filter is chosen as large as the number of measurements. As explained in Section II-C, Cadzow denoising exhibits a similar behaviour. This similarity is not fortuitous: both algorithms leverage a similar rank constraint which becomes more and more *selective* as the oversampling parameter increases. In contrast, GenFRI and LS-Cadzow are negatively affected by large oversampling parameters, due to numerical stability issues. For the critical bandwidth $2M + 1 = L$, CPGD notably outperforms GenFRI and LS-Cadzow by one to three orders of magnitude, and this even for PSNRs as low as $-20$ dB.

### B. Computational Complexity

As explained in Section III-A GenFRI has an overall computational complexity of $\mathcal{O}(N^3)$. For CPGD, the computational cost of each iteration is dominated by the successive projections onto $\mathcal{H}_K$ in the approximate proximal step, which are computed via an SVD –see Algorithm 1 and (14). At a cursory glance, it may seem that the overall complexity of CPGD is somewhat comparable to the one of GenFRI, since computing the SVD of a matrix with size $(N - P) \times (P + 1)$ has in general a computational complexity of $\mathcal{O}((N - P)^2(P + 1) + (P + 1)^3)$ [43] which reduces to $\mathcal{O}(N^3)$ when $P = M$. In practice however, projecting onto $\mathcal{H}_K$ does not require to perform a complete SVD since only the $K$ strongest eigenvalues and their associated eigenvectors are needed. This truncated SVD can be performed very efficiently by means of the *implicitly restarted Arnoldi method (IRAM)* [44], or the *implicitly restarted Lanczos method* for Hermitian matrices. When $K \ll P+1$, such methods are obviously much more efficient than a wasteful standard SVD. IRAM is moreover a *matrix-free* method [45]: it does not need the processed matrix to be stored in memory but simply requires an algorithm for performing matrix/vector products. In our context, since the truncated SVDs are exclusively performed on Toeplitz matrices, such matrix/vector products can be efficiently implemented by means of FFTs thanks to the convenient links between Toeplitz matrices and convolutions outlined in Section II-A.

[18]The performances of GenFRI were not investigated in this case, since the latter requires the forward matrix $\boldsymbol{G}$ to be injective which cannot be the case when $2M + 1 > L$.

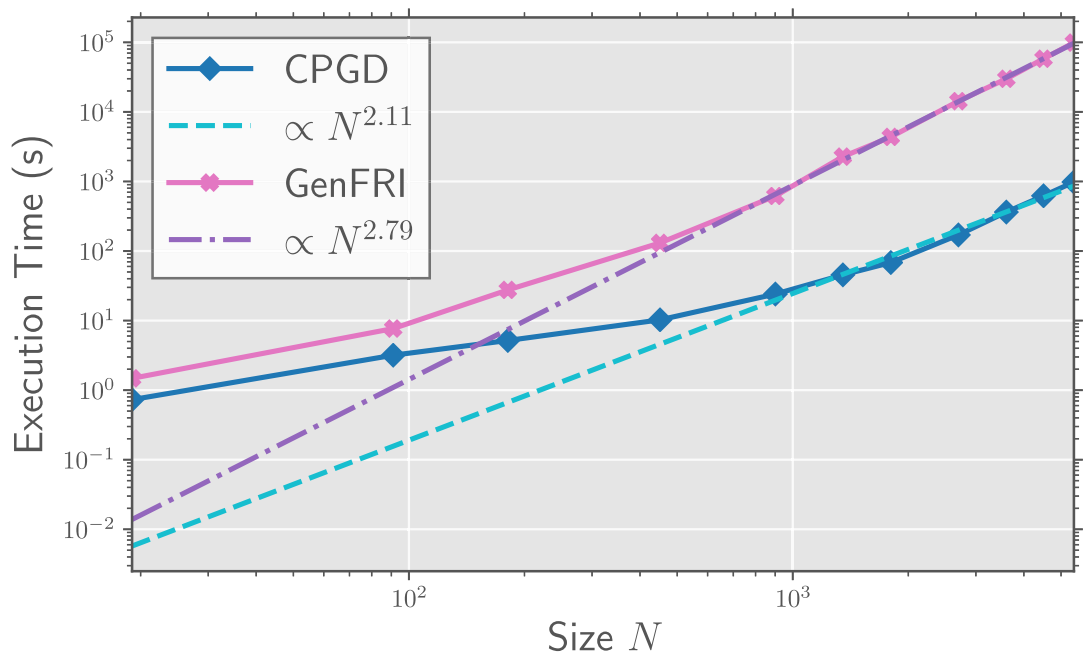

Fig. 5. Reconstruction times for CPGD and GenFRI for various bandwidth sizes $N$ and $K = 9$. The reported times are for a MacBook Pro (16-inch, 2019), Intel Core i7 (6 C/12 T) @ 2.6 GHz with 32 GB RAM. These results can be reproduced using the Python script `reproduce_execution_times.py` located in the directory `./benchmarking/` of our GitHub repository [30].

The CPGD implementation provided in our Github repository [30] leverages all these computational tricks. In Fig. 5, we show that our implementation of CPGD is considerably faster than GenFRI. Fig. 5 reports the reconstruction times of CPGD and GenFRI for $K = 9$ and bandwidth $N = 2\gamma M + 1 = L$ with $\gamma$ ranging from 1 to 300. To save computational time, the reconstruction times were scaled from the execution time of a single iteration of CPGD and GenFRI, assuming a typical number of iterations of 100 and 15×50=750 respectively. We observe that CPGD is always faster than GenFRI, sometimes by two orders of magnitude. Moreover, regressions performed in log-log scale reveal that CPGD scales as $N^{2.11}$ while GenFRI scales as $N^{2.79}$. Note that this difference in scaling behaviour is overlooked by the complexity analysis above.

## VI. Conclusion

We propose an implicit version of the generalised finite-rate-of-innovation (genFRI) problem for the recovery of the Fourier series coefficients of sparse Dirac streams with arbitrary linear sensing. This formulation relies on a novel regularisation term which enforces the annihilation of the recovered Fourier series coefficients without explicitly involving the unknown annihilating filter. The resulting non convex optimisation problem is consequently simpler and linear in the data. To solve it, we suggest a proximal gradient descent (PGD) algorithm which we prove converges towards a critical point of the objective function. We further introduce an inexact PGD method, coined *Cadzow plug-and-play gradient descent (CPGD)*, where the intractable proximal steps involved in PGD are approximated by means of alternating projections, akin to the popular Cadzow denoising algorithm. We outline the resemblance of CPGD to PnP methods used in image processing and prove its local fixed-point convergence under relatively weak assumptions. Considering the traditional irregular time sampling testbed, we demonstrate empirically that CPGD outperforms by several orders of magnitude the state-of-the-art GenFRI algorithm, both in terms of accuracy and reconstruction time.

For future work, we plan on investigating acceleration techniques for CPGD, such as approximate sketching-based eigenvalue decomposition methods [46], [47], more computationally efficient for large-scale problems. Applications of CPGD to acoustics and radio astronomy will also be investigated.

## Appendix A
## Proof of Proposition 1

Recall the definition of the proximal set associated to a point $\boldsymbol{x} \in \mathbb{C}^N$:

$$\operatorname{prox}_{\tau H}(\boldsymbol{x}) = \arg\min_{z\in\mathbb{C}^N}\left\{\frac{1}{2\tau}\|\boldsymbol{x}-z\|_2^2 + \iota_{\mathcal{H}_K}\left(T_P(z)\right) + \iota_{\mathbb{B}_\rho}(z)\right\}. \quad \text{(A.1)}$$

When mapped via the Toeplitzification operator $T_P$, the proximal set (A.1) becomes

$$\begin{aligned} T_P\left(\operatorname{prox}_{\tau H}(\boldsymbol{x})\right) &= \\ &= \left\{T_P(\check{x}),\ \check{x} \in \operatorname{prox}_{\tau H}(\boldsymbol{x})\right\} \\ &= \left\{\check{X} \in \mathbb{T}_P,\ T_P^\dagger(\check{X}) \in \operatorname{prox}_{\tau H}(\boldsymbol{x})\right\} \\ &= \arg\min_{\boldsymbol{Z}\in\mathbb{T}_P}\left\{\frac{1}{2\tau}\|T_P^\dagger(\boldsymbol{Z})-\boldsymbol{x}\|_2^2 + \iota_{\mathcal{H}_K}(\boldsymbol{Z}) + \iota_{\mathbb{B}_\rho}\left(T_P^\dagger(\boldsymbol{Z})\right)\right\} \\ &= \arg\min_{\boldsymbol{Z}\in\mathbb{T}_P\cap\mathcal{H}_K}\left\{\frac{1}{2\tau}\|T_P^\dagger(\boldsymbol{Z})-\boldsymbol{x}\|_2^2 + \iota_{\mathbb{B}_\rho}\left(T_P^\dagger(\boldsymbol{Z})\right)\right\}, \end{aligned} \quad \text{(A.2)}$$

where we have used the fact that $T_P^\dagger T_P(z) = z$ for all $z \in \mathbb{C}^N$. Define $\boldsymbol{W} = T_P(\operatorname{diag}(\Gamma^{-1/2})) \in \mathbb{T}_P$ where $\operatorname{diag}: \mathbb{C}^{N\times N} \to \mathbb{C}^N$ is the linear operator mapping a matrix onto its diagonal and $\Gamma = T_P^* T_P \in \mathbb{C}^{N\times N}$ is the diagonal and positive definite matrix given by (5). Then, the following relationships hold:

$$\Gamma^{-1/2}T_P^*(\boldsymbol{Z}) = T_P^*(\boldsymbol{W}\odot\boldsymbol{Z}), \quad \forall \boldsymbol{Z}\in\mathbb{C}^{(N-P)\times(P+1)}, \quad \text{(A.3)}$$

and

$$T_P(\Gamma^{-1/2}\boldsymbol{x}) = \boldsymbol{W}\odot T_P(\boldsymbol{x}), \quad \forall \boldsymbol{x}\in\mathbb{C}^N, \quad \text{(A.4)}$$

where $\odot$ denotes the *Hadamard product* for matrices. Using (A.3) and (A.4), we get moreover, $\forall \boldsymbol{Z} \in \mathbb{T}_P$:

$$\begin{aligned} \|T_P^\dagger(\boldsymbol{Z})-\boldsymbol{x}\|_2^2 &= \langle T_P^\dagger(\boldsymbol{Z})-\boldsymbol{x}, T_P^\dagger(\boldsymbol{Z})-\boldsymbol{x}\rangle_2 \\ &= \langle T_P^\dagger(\boldsymbol{Z}), T_P^\dagger(\boldsymbol{Z})\rangle_2 + \langle \boldsymbol{x},\boldsymbol{x}\rangle_2 \\ &\quad - \langle T_P^\dagger(\boldsymbol{Z}),\boldsymbol{x}\rangle_2 - \langle \boldsymbol{x}, T_P^\dagger(\boldsymbol{Z})\rangle_2 \\ &= \langle \Gamma^{-1}T_P^*(\boldsymbol{Z}), \Gamma^{-1}T_P^*(\boldsymbol{Z})\rangle_2 \\ &\quad + \langle \Gamma^{1/2}\Gamma^{-1/2}\boldsymbol{x}, \Gamma^{1/2}\Gamma^{-1/2}\boldsymbol{x}\rangle_2 \\ &\quad - \langle \Gamma^{-1}T_P^*(\boldsymbol{Z}),\boldsymbol{x}\rangle_2 - \langle \boldsymbol{x}, \Gamma^{-1}T_P^*(\boldsymbol{Z})\rangle_2 \end{aligned}$$

$$
\begin{aligned}
&= \langle \Gamma^{-1/2} T_P^*(\boldsymbol{W} \odot \boldsymbol{Z}), \Gamma^{-1/2} T_P^*(\boldsymbol{W} \odot \boldsymbol{Z}) \rangle_F \\
&\quad + \langle T_P^* T_P \Gamma^{-1/2} \boldsymbol{x}, \Gamma^{-1/2} \boldsymbol{x} \rangle_2 \\
&\quad - \langle T_P^*(\boldsymbol{W} \odot \boldsymbol{Z}), \Gamma^{-1/2} \boldsymbol{x} \rangle_2 \\
&\quad - \langle \Gamma^{-1/2} \boldsymbol{x}, T_P^*(\boldsymbol{W} \odot \boldsymbol{Z}) \rangle_2 \\
&= \langle \boldsymbol{W} \odot \boldsymbol{Z}, T_P \Gamma^{-1} T_P^*(\boldsymbol{W} \odot \boldsymbol{Z}) \rangle_F \\
&\quad + \langle T_P(\Gamma^{-1/2} \boldsymbol{x}), T_P(\Gamma^{-1/2} \boldsymbol{x}) \rangle_2 \\
&\quad - \langle \boldsymbol{W} \odot \boldsymbol{Z}, T_P(\Gamma^{-1/2} \boldsymbol{x}) \rangle_F \\
&\quad - \langle T_P(\Gamma^{-1/2} \boldsymbol{x}), \boldsymbol{W} \odot \boldsymbol{Z} \rangle_F \\
&= \langle \boldsymbol{W} \odot \boldsymbol{Z}, \boldsymbol{W} \odot \boldsymbol{Z} \rangle_F \\
&\quad + \langle \boldsymbol{W} \odot T_P(\boldsymbol{x}), \boldsymbol{W} \odot T_P(\boldsymbol{x}) \rangle_F \\
&\quad - \langle \boldsymbol{W} \odot \boldsymbol{Z}, \boldsymbol{W} \odot T_P(\boldsymbol{x}) \rangle_F \\
&\quad - \langle \boldsymbol{W} \odot T_P(\boldsymbol{x}), \boldsymbol{W} \odot \boldsymbol{Z} \rangle_F \\
&= \|\boldsymbol{W} \odot \boldsymbol{Z}\|_F^2 + \|\boldsymbol{W} \odot T_P(\boldsymbol{x})\|_F^2 \\
&\quad - 2\Re \left( \langle \boldsymbol{W} \odot \boldsymbol{Z}, \boldsymbol{W} \odot T_P(\boldsymbol{x}) \rangle_F \right) \\
&= \|\boldsymbol{W} \odot (\boldsymbol{Z} - T_P(\boldsymbol{x}))\|_F^2, \qquad \text{(A.5)}
\end{aligned}
$$

where we have used the fact that $\Gamma = \Gamma^H = T_P^* T_P$, $\boldsymbol{W} \odot \boldsymbol{Z} \in \mathbb{T}_P$ for every $\boldsymbol{Z} \in \mathbb{T}_P$ and $T_P \Gamma^{-1} T_P^* = \Pi_{\mathbb{T}_P}$ (see Appendices E and F). With similar arguments, we have $\forall \boldsymbol{Z} \in \mathbb{T}_P$:

$$
\begin{aligned}
\left\| T_P^\dagger(\boldsymbol{Z}) \right\|_2 \le \rho &\Leftrightarrow \sqrt{\langle T_P^\dagger(\boldsymbol{Z}), T_P^\dagger(\boldsymbol{Z}) \rangle_2} \le \rho \\
&\Leftrightarrow \sqrt{\langle \boldsymbol{W} \odot \boldsymbol{Z}, \boldsymbol{W} \odot \boldsymbol{Z} \rangle_F} \le \rho \\
&\Leftrightarrow \|\boldsymbol{W} \odot \boldsymbol{Z}\|_F \le \rho,
\end{aligned}
$$

so that

$$
\iota_{\mathbb{B}_\rho}\left( T_P^\dagger(\boldsymbol{Z}) \right) = \iota_{\mathbb{B}_\rho^{\boldsymbol{W}}}(\boldsymbol{Z}), \quad \forall \boldsymbol{Z} \in \mathbb{T}_P, \qquad \text{(A.6)}
$$

where $\mathbb{B}_\rho^{\boldsymbol{W}} := \{\boldsymbol{Z} \in \mathbb{C}^{(N-P)\times(P+1)} : \|\boldsymbol{W} \odot \boldsymbol{Z}\|_F \le \rho\}$.

Plugging (A.5) and (A.6) into (A.2) hence yields

$$
\begin{aligned}
&T_P\left(\mathrm{prox}_{\tau H}(\boldsymbol{x})\right) = \\
&= \arg \min_{\boldsymbol{Z} \in \mathbb{T}_P \cap \mathcal{H}_K} \left\{ \frac{1}{2\tau} \|T_P^\dagger(\boldsymbol{Z}) - \boldsymbol{x}\|_2^2 + \iota_{\mathbb{B}_\rho}\left( T_P^\dagger(\boldsymbol{Z}) \right) \right\} \\
&= \arg \min_{\boldsymbol{Z} \in \mathbb{T}_P \cap \mathcal{H}_K} \left\{ \frac{1}{2\tau} \|\boldsymbol{W} \odot (\boldsymbol{Z} - T_P(\boldsymbol{x}))\|_F^2 + \iota_{\mathbb{B}_\rho^{\boldsymbol{W}}}(\boldsymbol{Z}) \right\} \\
&= \arg \min_{\boldsymbol{Z} \in \mathbb{T}_P \cap \mathcal{H}_K \cap \mathbb{B}_\rho^{\boldsymbol{W}}} \left\{ \frac{1}{2\tau} \|\boldsymbol{W} \odot (\boldsymbol{Z} - T_P(\boldsymbol{x}))\|_F^2 \right\} \\
&= \arg \min_{\boldsymbol{Z} \in \mathbb{T}_P \cap \mathcal{H}_K \cap \mathbb{B}_\rho^{\boldsymbol{W}}} \|\boldsymbol{W} \odot (\boldsymbol{Z} - T_P(\boldsymbol{x}))\|_F \\
&= \Pi^{\boldsymbol{W}}_{\mathbb{T}_P \cap \mathcal{H}_K \cap \mathbb{B}_\rho^{\boldsymbol{W}}} T_P(\boldsymbol{x}). \qquad \text{(A.7)}
\end{aligned}
$$

Using the fact that $T_P^\dagger T_P = \boldsymbol{I}_N$ we can finally rewrite (A.7) as

$$
\mathrm{prox}_{\tau H}(\boldsymbol{x}) = T_P^\dagger \Pi^{\boldsymbol{W}}_{\mathbb{T}_P \cap \mathcal{H}_K \cap \mathbb{B}_\rho^{\boldsymbol{W}}} T_P(\boldsymbol{x}),
$$

which completes the proof.

## Appendix B
## Proof of Theorem 3

The proof of Theorem 3 relies on the four lemmas hereafter. The first lemma shows that gradient descent is Lipschitz continuous, and exhibits step size ranges for which it is also a *$\mu$-contraction* –i.e., the Lipschitz constant is strictly smaller than $1/\sqrt{\mu+1}$ for some $\mu \ge 0$. This is a generalisation of a famous result in optimisation [48], [49].

*Lemma B.1 ($\mu$-Contractive Gradient Descent):* Let $\boldsymbol{G} \in \mathbb{C}^{L\times N}$ be injective, and define

$$
\alpha := 2\lambda_{\min}\left(\boldsymbol{G}^H \boldsymbol{G}\right), \qquad \text{(B.1)}
$$

$$
\beta := 2\lambda_{\max}\left(\boldsymbol{G}^H \boldsymbol{G}\right), \qquad \text{(B.2)}
$$

where $\lambda_{\min}(\boldsymbol{M})$ and $\lambda_{\max}(\boldsymbol{M})$ denote the minimum and maximum eigenvalue of a matrix $\boldsymbol{M}$ respectively. Let $\tau \in \mathbb{R}_+$ be a positive constant and consider the linear map

$$
D_\tau : \begin{cases} \mathbb{C}^N \to \mathbb{C}^N, \\ \boldsymbol{x} \mapsto \boldsymbol{x} - 2\tau \boldsymbol{G}^H (\boldsymbol{G}\boldsymbol{x} - \boldsymbol{y}), \end{cases} \qquad \text{(B.3)}
$$

for some $\boldsymbol{y} \in \mathbb{C}^L$. Then, $D_\tau$ is Lipschitz continuous:

$$
\|D_\tau(\boldsymbol{x}) - D_\tau(\boldsymbol{z})\|_2 \le L_\tau \|\boldsymbol{x} - \boldsymbol{z}\|_2, \qquad \forall \boldsymbol{x}, \boldsymbol{z} \in \mathbb{C}^N,
$$

with Lipschitz contant:

$$
L_\tau = \max\{|1 - \tau\alpha|, |1 - \tau\beta|\}. \qquad \text{(B.4)}
$$

For $\mu \ge 0$ moreover, $D_\tau$ is $\mu$-contractive, i.e., $0 \le L_\tau < 1/\sqrt{\mu+1}$, for

$$
\frac{1}{\beta}\left(1 - \frac{1}{\sqrt{\mu+1}}\right) < \tau < \frac{1}{\beta}\left(1 + \frac{1}{\sqrt{\mu+1}}\right). \qquad \text{(B.5)}
$$

*Proof:* We have

$$
\begin{aligned}
\|D_\tau(\boldsymbol{x}) - D_\tau(\boldsymbol{z})\|_2 &= \left\| (\boldsymbol{I}_N - 2\tau \boldsymbol{G}^H \boldsymbol{G})(\boldsymbol{x} - \boldsymbol{z}) \right\|_2 \\
&\le \left\| \boldsymbol{I}_N - 2\tau \boldsymbol{G}^H \boldsymbol{G} \right\|_2 \|\boldsymbol{x} - \boldsymbol{z}\|_2 \\
&= L_\tau \|\boldsymbol{x} - \boldsymbol{z}\|_2,
\end{aligned}
$$

where the Lipschitz constant $L_\tau := \|\boldsymbol{I}_N - 2\tau \boldsymbol{G}^H \boldsymbol{G}\|_2 > 0$ is the spectral norm of $\boldsymbol{I}_N - 2\tau \boldsymbol{G}^H \boldsymbol{G}$. Note that since $\boldsymbol{G}$ is injective, $\boldsymbol{G}^H \boldsymbol{G}$ is positive definite and hence we easily get [49] that the eigenvalues of $\boldsymbol{I}_N - 2\tau \boldsymbol{G}^H \boldsymbol{G}$ are contained in the interval $[1 - \tau\beta, 1 - \tau\alpha]$, where $\beta \ge \alpha > 0$ are defined in (B.1) and (B.2) respectively. Its spectral norm is hence given by:

$$
\left\| \boldsymbol{I}_N - 2\tau \boldsymbol{G}^H \boldsymbol{G} \right\|_2 = \max\{|1 - \tau\alpha|, |1 - \tau\beta|\},
$$

which proves (B.4).

For $\tau$ verifying (B.5), we have moreover:

$$
\begin{aligned}
&1 - \frac{1}{\sqrt{\mu+1}} < \tau\beta < 1 + \frac{1}{\sqrt{\mu+1}} \\
&\Leftrightarrow -\frac{1}{\sqrt{\mu+1}} < 1 - \tau\beta < \frac{1}{\sqrt{\mu+1}} \\
&\Leftrightarrow |1 - \tau\beta| < \frac{1}{\sqrt{\mu+1}}
\end{aligned}
$$

and similarly since $\alpha \leq \beta$:

$$\frac{\alpha}{\beta}\left(1 - \frac{1}{\sqrt{\mu+1}}\right) < \tau\alpha < \frac{\alpha}{\beta}\left(1 + \frac{1}{\sqrt{\mu+1}}\right)$$

$$\Leftrightarrow 1 - \frac{1}{\sqrt{\mu+1}} < \tau\alpha < 1 + \frac{1}{\sqrt{\mu+1}}$$

$$\Leftrightarrow |1 - \tau\alpha| < \frac{1}{\sqrt{\mu+1}},$$

which shows the $\mu$-contractivity of $D_\tau$. ■

The second lemma states that in a Hilbert space, projection maps onto closed convex sets are non-expansive. This is a known result from approximation theory [50], [51].

*Lemma B.2 (Non-Expansiveness of Closed Convex Projections):* Let $\mathcal{H}$ be some Hilbert space with some inner-product norm $\|\cdot\|$ and $\mathcal{C} \subset \mathcal{H}$ a *closed*, *convex* set. Then the projection map onto $\mathcal{C}$, defined as

$$\Pi_{\mathcal{C}}(x) = \arg\min_{z \in \mathcal{C}} \|x - z\|, \qquad \forall x \in \mathcal{H},$$

is *non-expansive*, i.e.,

$$\|\Pi_{\mathcal{C}}(x) - \Pi_{\mathcal{C}}(z)\| \leq \|x - z\|, \quad \forall x, z \in \mathcal{H}.$$

*Proof:* Lemma B.2 is proven in [50, Theorem 5.5]. ■

The third lemma states that the singular value projection map $\Pi_{\mathcal{H}_k}$ is locally non-expansive in every neighbourhood of the manifold of matrices with rank exactly $k$.

*Lemma B.3 (Local Non-Expansiveness of the Singular Value Projection):* Let $\mathbb{C}^{m\times n}$ be the space of complex-valued rectangular matrices of size $m \times n$, and $\mathcal{H}_k \subset \mathbb{C}^{m\times n}$, $\mathcal{R}_k \subset \mathbb{C}^{m\times n}$ the sets of matrices with rank at most and exactly $k \leq \max\{m, n\}$ respectively. Denote further by $\Pi_{\mathcal{H}_k}$ the projection map onto $\mathcal{H}_k$ given in (14). Then, for every $\boldsymbol{R} \in \mathcal{R}_k$, the map $\Pi_{\mathcal{H}_k}$ is well-defined (single-valued) and *locally non-expansive*

$$\|\Pi_{\mathcal{H}_k}(\boldsymbol{X}) - \Pi_{\mathcal{H}_k}(\boldsymbol{Z})\|_F \leq \|\boldsymbol{X} - \boldsymbol{Z}\|_F, \quad \forall \boldsymbol{X}, \boldsymbol{Y} \in \mathcal{U},$$

for some neighbourhood[19] $\mathcal{U} \ni \boldsymbol{R}$.

*Proof:* Since $\mathcal{R}_k$ is dense in $\mathcal{H}_k$ [34, Proposition 2.1], we have $\Pi_{\mathcal{H}_k} = \Pi_{\mathcal{R}_k}$ in a neighbourhood $\mathcal{W}$ of every $\boldsymbol{R} \in \mathcal{R}_k$ (see [52, Example 2.3] for a detailed proof of this fact). Moreover, [53, Lemma 3] tells us that, for every $\boldsymbol{R} \in \mathcal{R}_k$, $\Pi_{\mathcal{R}_k}$ is, in a neighbourhood $\mathcal{U} \ni \boldsymbol{R}$ such that $\mathcal{U} \subset \mathcal{W}$, well-defined (single-valued), continuous and differentiable, with gradient given by: $\nabla\Pi_{\mathcal{R}_k} = \Pi_{\mathbb{T}_{\mathcal{R}_k}(\boldsymbol{R})}$ where $\mathbb{T}_{\mathcal{R}_k}(\boldsymbol{R}) \subset \mathbb{C}^{m\times n}$ is the tangent plane of the manifold $\mathcal{R}_k$ in $\boldsymbol{R}$ (see [52, Example 2.2]). Since $\mathbb{T}_{\mathcal{R}_k}(\boldsymbol{R})$ is by definition a linear subspace of $\mathbb{C}^{m\times n}$, the orthogonal projection operator $\Pi_{\mathbb{T}_{\mathcal{R}_k}(\boldsymbol{R})}$ is bounded with unit spectral norm. The map $\Pi_{\mathcal{R}_k} = \Pi_{\mathcal{H}_k}$ is consequently 1-Lipschitz continuous (i.e., non-expansive) with respect to the Frobenius norm in the neighbourhood $\mathcal{U}$ of $\boldsymbol{R} \in \mathcal{R}_k$. ■

The last lemma finally, makes use of Lemmas B.2 and B.3 to show that the denoising operator $H_n(\boldsymbol{x}) = T_P^\dagger[\Pi_{\mathbb{T}_P}\Pi_{\mathcal{H}_K}\Pi_{\mathbb{B}_\rho^{\boldsymbol{W}}}^{\boldsymbol{W}}]^n T_P(\boldsymbol{x})$ is locally Lipschitz continuous:

*Lemma B.4 (Local Lipschitz Continuity of Denoiser):* Let $\mathbb{C}^{(N-P)\times(P+1)}$ be the space of complex-valued rectangular matrices of size $(N-P) \times (P+1)$, $P \leq \lfloor N/2 \rfloor$, and $\mathcal{H}_K \subset \mathbb{C}^{(N-P)\times(P+1)}$, $\mathcal{R}_K \subset \mathbb{C}^{(N-P)\times(P+1)}$ the sets of matrices with rank at most and exactly $K \leq P$ respectively. Let

$$H_n(\boldsymbol{x}) := T_P^\dagger\left[\Pi_{\mathbb{T}_P}\Pi_{\mathcal{H}_K}\Pi_{\mathbb{B}_\rho^{\boldsymbol{W}}}^{\boldsymbol{W}}\right]^n T_P(\boldsymbol{x}), \quad \forall \boldsymbol{x} \in \mathbb{C}^N,$$

be the approximate proximal operator (33). Then, $H_n$ is locally well-defined (single-valued) and $\sqrt{P+1}$-Lipschitz continuous

$$\|H_n(\boldsymbol{x}) - H_n(\boldsymbol{z})\|_2 \leq \sqrt{P+1}\|\boldsymbol{x} - \boldsymbol{z}\|_2,$$

for all $\boldsymbol{x}, \boldsymbol{z} \in \mathbb{C}^N$ such that $T_P(\boldsymbol{x})$, $T_P(\boldsymbol{z})$ are in some neighbourhood of some matrix $\boldsymbol{R} \in \mathcal{R}_K$.

*Proof:* First, we have, for all $\boldsymbol{x}, \boldsymbol{z} \in \mathbb{C}^N$:

$$\|H_n(\boldsymbol{x}) - H_n(\boldsymbol{z})\|_2 = \left\|T_P^\dagger(\Pi_n T_P(\boldsymbol{x}) - \Pi_n T_P(\boldsymbol{z}))\right\|_2, \tag{B.6}$$

where $\Pi_n = [\Pi_{\mathbb{T}_P}\Pi_{\mathcal{H}_K}\Pi_{\mathbb{B}_\rho^{\boldsymbol{W}}}^{\boldsymbol{W}}]^n$. As shown in (A.5) we have moreover, for $\boldsymbol{X} \in \mathbb{T}_P$,

$$\left\|T_P^\dagger(\boldsymbol{X})\right\|_2^2 = \|\boldsymbol{W} \odot \boldsymbol{X}\|_F^2 = \sum_{i=1}^{N-P}\sum_{j=1}^{P+1} W_{ij}^2 |X_{ij}|^2$$
$$\leq \left\|\boldsymbol{W}^{\odot 2}\right\|_\infty \|\boldsymbol{X}\|_F^2,$$

where $\boldsymbol{W} = T_P(\mathrm{diag}(\Gamma^{-1/2})) \in \mathbb{T}_P$ is as in Proposition 1. From the definition of $\Gamma$ in (5), it is moreover easy to show that $\|\boldsymbol{W}^{\odot 2}\|_\infty = \|\Gamma^{-1}\|_\infty = 1$ and hence

$$\left\|T_P^\dagger(\boldsymbol{X})\right\|_2^2 \leq \|\boldsymbol{X}\|_F^2, \qquad \forall \boldsymbol{X} \in \mathbb{T}_P.$$

Since the range of $\Pi_n$ is $\mathbb{T}_P$, (B.6) becomes:

$$\|H_n(\boldsymbol{x}) - H_n(\boldsymbol{z})\|_2 \leq \|\Pi_n T_P(\boldsymbol{x}) - \Pi_n T_P(\boldsymbol{z})\|_F.$$

Assuming now that $T_P(\boldsymbol{x})$ and $T_P(\boldsymbol{z})$ are in some neighbourhood of some point $\boldsymbol{R} \in \mathcal{R}_K$, we can invoke Lemmas B.2 and B.3 recursively to obtain:

$$\begin{aligned}\|\Pi_n T_P(\boldsymbol{x}) - \Pi_n T_P(\boldsymbol{z})\|_F &\leq \|T_P(\boldsymbol{x}) - T_P(\boldsymbol{z})\|_F \\ &= \left\|\Gamma^{1/2}(\boldsymbol{x} - \boldsymbol{z})\right\|_2 \\ &\leq \left\|\Gamma^{1/2}\right\|_2 \|\boldsymbol{x} - \boldsymbol{z}\|_2,\end{aligned}$$

where we have used: $\|T_P(\boldsymbol{x})\|_F^2 = \langle T_P(\boldsymbol{x}), T_P(\boldsymbol{x})\rangle_F = \langle T_P^* T_P(\boldsymbol{x}), \boldsymbol{x}\rangle_2 = \|\Gamma^{1/2}\boldsymbol{x}\|_2$, $\forall \boldsymbol{x} \in \mathbb{C}^N$. From the definition of $\Gamma$ in (5), we can easily show that $\|\Gamma^{1/2}\|_2 = \sqrt{P+1}$ which finally yields

$$\|H_n(\boldsymbol{x}) - H_n(\boldsymbol{z})\|_2 \leq \sqrt{P+1}\,\|\boldsymbol{x} - \boldsymbol{z}\|_2,$$

for all $\boldsymbol{x}, \boldsymbol{z} \in \mathbb{C}^N$ such that $T_P(\boldsymbol{x})$, $T_P(\boldsymbol{z})$ are in some neighbourhood of some matrix $\boldsymbol{R} \in \mathcal{R}_K$. ■

We are now ready to show Theorem 3. Let

$$U_{\tau,n}(\boldsymbol{x}) := H_n(\boldsymbol{x} - \tau\nabla F(\boldsymbol{x})) = H_n(D_\tau(\boldsymbol{x})), \qquad \boldsymbol{x} \in \mathbb{C}^N.$$

Then, for every $\boldsymbol{x}, \boldsymbol{z} \in \mathbb{C}^N$ such that $T_P(\boldsymbol{x})$, $T_P(\boldsymbol{z})$ are in some neighbourhood of some matrix $\boldsymbol{R} \in \mathcal{R}_K$, $U_{\tau,n}$ is locally Lipschitz continuous as composition between two (locally) Lipschitz continuous functions $H_n$ and $D_\tau$, see Lemmas B.4 and B.1 respectively. Moreover, the Lipschitz constant is the product of the Lipschitz constants of $H_n$ and $D_\tau$, $\sqrt{P+1}$ and $L_\tau$ in (B.4) respectively. We have therefore

$$\left\|U_{\tau,n}(\boldsymbol{x}) - U_{\tau,n}(\boldsymbol{z})\right\|_2 \le \sqrt{P+1} L_\tau \|\boldsymbol{x}-\boldsymbol{z}\|_2,$$

for all $\boldsymbol{x}, \boldsymbol{z} \in \mathbb{C}^N$ such that $T_P(\boldsymbol{x}), T_P(\boldsymbol{z})$ are in some neighbourhood of some matrix $\boldsymbol{R} \in \mathcal{R}_K$. Finally, for

$$\frac{1}{\beta}\left(1 - \frac{1}{\sqrt{P+1}}\right) < \tau < \frac{1}{\beta}\left(1 + \frac{1}{\sqrt{P+1}}\right)$$

we have from Lemma B.1 that $L_\tau < 1/\sqrt{P+1}$ and hence finally

$$\left\|U_{\tau,n}(\boldsymbol{x}) - U_{\tau,n}(\boldsymbol{z})\right\|_2 < \|\boldsymbol{x}-\boldsymbol{z}\|_2,$$

which shows than $U_{\tau,n}$ is locally contractive. ∎

[19] If $X$ is a topological space and $p$ is a point in $X$, a *neighbourhood* of $p$ is a subset $V$ of $X$ that includes an open set $U$ containing $p$, $p \in U \subseteq V$.

## Appendix C
## Proof of Corollary 1

First, we note that from Theorem 3, we have under the assumptions of the corollary that

$$\|\boldsymbol{x}_{k+1} - \boldsymbol{x}_k\|_2 = \|U_{\tau,n}(\boldsymbol{x}_k) - U_{\tau,n}(\boldsymbol{x}_{k-1})\|_2 \le \widetilde{L}_\tau \|\boldsymbol{x}_k - \boldsymbol{x}_{k-1}\|_2,$$

for all $k \ge 1$ and hence by induction

$$\|\boldsymbol{x}_{k+1} - \boldsymbol{x}_k\|_2 \le \widetilde{L}_\tau^k \|\boldsymbol{x}_1 - \boldsymbol{x}_0\|_2, \quad \forall k \ge 1. \tag{C.1}$$

Since $\tau$ is assumed to satisfy (36) we have moreover $0 < \widetilde{L}_\tau < 1$. We deduce hence from (C.1) that $\{\boldsymbol{x}_k\}_{k\in\mathbb{N}}$ is a Cauchy sequence. Let $j, k \in \mathbb{N}$ with $j > k$:

$$\begin{aligned}
\|\boldsymbol{x}_j - \boldsymbol{x}_k\|_2 &\le \sum_{m=k}^{j-1} \|\boldsymbol{x}_{m+1} - \boldsymbol{x}_m\|_2 \\
&\le \sum_{m=k}^{j-1} \widetilde{L}_\tau^m \|\boldsymbol{x}_1 - \boldsymbol{x}_0\|_2 \\
&= \|\boldsymbol{x}_1 - \boldsymbol{x}_0\|_2 \widetilde{L}_\tau^k \sum_{m=0}^{j-1-k} \widetilde{L}_\tau^m \\
&\le \|\boldsymbol{x}_1 - \boldsymbol{x}_0\|_2 \widetilde{L}_\tau^k \sum_{m=0}^{\infty} \widetilde{L}_\tau^m \\
&= \frac{\widetilde{L}_\tau^k}{1-\widetilde{L}_\tau} \|\boldsymbol{x}_1 - \boldsymbol{x}_0\|_2.
\end{aligned} \tag{C.2}$$

For every $\epsilon > 0$, we can choose a $J \in \mathbb{N}$ such that

$$\widetilde{L}_\tau^J < \frac{\epsilon(1-\widetilde{L}_\tau)}{\|\boldsymbol{x}_1 - \boldsymbol{x}_0\|_2},$$

and hence for all $j > k > J$

$$\|\boldsymbol{x}_j - \boldsymbol{x}_k\|_2 < \epsilon.$$

The sequence $\{\boldsymbol{x}_k\}_{k\in\mathbb{N}}$ is hence a Cauchy sequence, and since $\mathbb{C}^N$ is complete, it converges towards a limit point $\boldsymbol{x}_\star \in \mathbb{C}^N$. We have moreover, since $U_{\tau,n}$ is continuous

$$\begin{aligned}
\boldsymbol{x}_\star &= \lim_{k\to\infty} \boldsymbol{x}_k = \lim_{k\to\infty} U_{\tau,n}(\boldsymbol{x}_{k-1}) = U_{\tau,n}\left(\lim_{k\to\infty} \boldsymbol{x}_{k-1}\right) \\
&= U_{\tau,n}(\boldsymbol{x}_\star),
\end{aligned}$$

and hence $\boldsymbol{x}_\star$ is a fixed-point of $U_{\tau,n}$. Note moreover that, from (C.2) we get

$$\begin{aligned}
\|\boldsymbol{x}_\star - \boldsymbol{x}_k\|_2 &= \lim_{j\to+\infty} \|\boldsymbol{x}_j - \boldsymbol{x}_k\|_2 \\
&\le \lim_{j\to+\infty} \|\boldsymbol{x}_1 - \boldsymbol{x}_0\|_2 \widetilde{L}_\tau^k \sum_{m=0}^{j-1-k} \widetilde{L}_\tau^m \\
&= \|\boldsymbol{x}_1 - \boldsymbol{x}_0\|_2 \widetilde{L}_\tau^k \sum_{m=0}^{+\infty} \widetilde{L}_\tau^m \\
&= \frac{\widetilde{L}_\tau^k}{1-\widetilde{L}_\tau} \|\boldsymbol{x}_1 - \boldsymbol{x}_0\|_2,
\end{aligned}$$

which proves (38) of Corollary 1. ∎

## Acknowledgment

The authors are grateful to Sepand Kashani, Julien Fageot, Adrian Jarret and the anonymous referees for their valuable comments and suggestions, which greatly helped improving the overall quality of the paper.